\documentclass[preprint,12pt]{elsarticle}

\usepackage{amssymb}
\usepackage{amsmath}
\usepackage{hyperref}

\usepackage[utf8]{inputenc}
\usepackage{
amsmath,
amssymb,
amsthm,
appendix,
bm,
bbm,
caption,
color,
epsf,
enumitem,
float,
graphicx,
hyperref,
lineno,
listings,
mathtools,
mathrsfs,
siunitx,
subcaption,
tikz,
titletoc,
url,
ulem,
xcolor
}

\usepackage[margin=1in,top=1in]{geometry}
\hypersetup{
    colorlinks=true, 
    linktoc=all,     
    linkcolor=black,  
}

\numberwithin{equation}{section}

\tikzset{every label/.style={font=\footnotesize,inner sep=1pt}}
\DeclareCaptionFormat{myformat}{#1#2#3\hrulefill}
\DeclareCaptionFormat{myformat2}{#1#2#3}
\newcommand{\Ebb}{\mathbb{E}}

\newcommand{\Lbb}{\mathbb{L}}

\newcommand{\Rbb}{\mathbb{R}}
\newcommand{\Sbb}{\mathbb{S}}

\newcommand{\Ubb}{\mathbb{U}}

\newcommand{\Abf}{\mathbf{A}}
\newcommand{\bbf}{\mathbf{b}}
\newcommand{\Bbf}{\mathbf{B}}

\newcommand{\Cbf}{\mathbf{C}}

\newcommand{\Dbf}{\mathbf{D}}
\newcommand{\ebf}{\mathbf{e}}

\newcommand{\Fbf}{\mathbf{F}}

\newcommand{\Gbf}{\mathbf{G}}

\newcommand{\Ibf}{\mathbf{I}}

\newcommand{\pbf}{\mathbf{p}}

\newcommand{\qbf}{\mathbf{q}}

\newcommand{\Sbf}{\mathbf{S}}
\newcommand{\tbf}{\mathbf{t}}
\newcommand{\Tbf}{\mathbf{T}}
\newcommand{\ubf}{\mathbf{u}}
\newcommand{\Ubf}{\mathbf{U}}

\newcommand{\wbf}{\mathbf{w}}
\newcommand{\Wbf}{\mathbf{W}}
\newcommand{\xbf}{\mathbf{x}}
\newcommand{\Xbf}{\mathbf{X}}

\newcommand{\ep}{\epsilon}

\newcommand{\CalE}{{\mathcal{E}}}
\newcommand{\CalF}{{\mathcal{F}}}

\newcommand{\CalI}{{\mathcal{I}}}

\newcommand{\CalL}{{\mathcal{L}}}

\newcommand{\CalO}{{\mathcal{O}}}
\newcommand{\CalQ}{{\mathcal{Q}}}

\newcommand{\CalS}{{\mathcal{S}}}

\newcommand{\bOmega}{\boldsymbol{\Omega}}
\renewcommand{\d}{\mathrm{d}}

\newcommand{\nrm}[1]{\left\Vert {#1} \right\Vert}

\newcommand{\supp}[1]{\text{supp}\left(#1\right)}

\renewcommand{\tilde}[1]{\widetilde{#1}}

\DeclareMathOperator*{\argmin}{argmin}

\newcommand{\lan}{\left\langle}
\newcommand{\ran}{\right\rangle}

\newcommand{\wstar}{\wbf^\star}
\newcommand{\what}{{\widehat{\wbf}}}

\newcommand{\bi}[1]{\textbf{\textit{#1}}}

\journal{J. Comp. Phys}

\begin{document}

\begin{frontmatter}



\title{Learning interpretable closures for thermal radiation transport in optically-thin media using WSINDy}


\author{Daniel Messenger, Ben Southworth, Hans Hammer, Luis Chacon} 

\affiliation{organization={Theoretical Division, Los Alamos National Laboratory},
            city={Los Alamos},
            postcode={87544}, 
            state={New Mexico},
            country={USA}}

\begin{abstract}
We introduce an equation learning framework to identify a closed set of equations for moment quantities in 1D thermal radiation transport (TRT) in optically thin media. While optically thick media admits a well-known diffusive closure, the utility of moment closures in providing accurate low-dimensional surrogates for TRT in optically thin media is unclear, as the mean-free path of photons is large and 
the radiation flux is far from its Fickean limit. 
Here, we demonstrate the viability of using weak-form equation learning to close the system of equations for the energy density, radiation flux, and temperature in optically thin TRT. We show that the WSINDy algorithm (Weak-form Sparse Identification of Nonlinear Dynamics), together with an advantageous change of variables and an auxiliary equation for the radiation-energy-weighted opacity, enables robust and efficient identification of closures that preserve many desired physical properties from the high fidelity system, including hyperbolicity, rotational symmetry, black-body equilibria, and linear stability of black-body equilibria, all of which manifest as library constraints or convex constraints on the closure coefficients. Crucially, the weak form enables closures to be learned from simulation data with ray effects and particle noise, which then do not appear in simulations of the resulting closed moment system. 
Finally, we demonstrate that our closure models can be extrapolated in the key system parameters of drive temperature $T_{in}$ and scalar opacity $\gamma$, and this extrapolation is to an extent quantifiable by a Knudsen-like dimensionless parameter.  

\end{abstract}


\begin{highlights}
\item Parametrized closures for optically thin thermal radiation transport
\item Weak formulation mitigates ray effects and particle noise in kinetic data
\item Hyperbolicity and key physical properties enforced as convex constraints
\item Extrapolation in parameter space and spatiotemporal learning domain
\end{highlights}


\end{frontmatter}



\section{Introduction}

    Even in well-studied physical systems with known governing equations (e.g.\ Newtonian molecular dynamics, turbulent fluids, plasma physics), myriad challenges arise in running simulations due to the curse of dimensionality, whereby the high fidelity system is impossible to resolve at engineering scales using state-of-the-art computing resources. Thermal radiation transport (TRT), a continuum model for the interactions between photons and matter, is one such physical theory
    that represents the primary computational bottleneck in many multiphysics simulations, including radiation hydrodynamics (rad-hydro, used to model inertial confinement fusion experiments \cite{slutz2006integrated,colvin2011role,pak2020impact,haines2021constraining,kim2023evaluation,kurzer2024radiation,marinak2024numerical}), massive star envelope formation \cite{jiang2023three}, and general astrophysical plasmas \cite{rochau2014zapp}. The key state variables in TRT are the specific intensity of radiation $I = I(\xbf,\bOmega,\nu,t)$, a field over space $\xbf\in \Rbb^3$, angle $\bOmega\in \Sbb^2$, frequency $\nu \in \Rbb_+$, and time $t$, and the material temperature $T = T(\xbf,t)$. The partial differential equations (PDEs) governing TRT are given below in \eqref{eq:ho}.
    
    Closure modeling is an ubiquitous approach for breaking the curse of dimensionality, whereby only averaged physical quantities are directly modeled. This tends to remove dimensions from phase space while increasing the state-space dimension, a much more favorable tradeoff considering the computational complexity of $\CalO(n N^d)$ for resolving $n$ state variables over a $d$-dimensional phase space with $N$ points per dimension.  Implementations of numerical rad-hydro solvers typically allow users to select from a suite of closure models for TRT depending on how diffusive, opaque, isotropic, or gray the radiation dynamics are. For gray isotropic TRT in optically-thick media, the system closes in the radiation energy density $E$ and temperature $T$, effectively a zeroth moment closure ($E$ is the zeroth moment of the radiation intensity $I$). This model can be extended to multifrequency scenarios, modeling multiple photon frequencies $\nu$, where energy groups $E_\nu$ all couple through the temperature equation. When the medium is not optically thick, but the radiation field is still approximately isotropic, a diffusive approximation cannot be justified, but a first-order moment closure involving the radiation flux $F$ (first moment of $I$) is appropriate. This leads to the commonly employed $P_1$ and $P_{1/3}$ models \cite{olson2000diffusion}. When the dynamics are optically {\it thin}, the radiation field is almost certainly anisotropic, such that a nontrivial closure for the second moment of $I$, known as the Eddington tensor $\CalE$, is required. When the physics falls outside of any known closure regime, rad-hydro simulations must solve the full TRT system, typically by employing particle methods for photon transport such as Implicit Monte Carlo \cite{fleck1971implicit,haines2020cross,marinak2024numerical} or Deterministic Particle (DP) methods \cite{park2019multigroup,hammer2019multi}, or  Eulerian methods such as finite elements  \cite{liu2006finite,lou2019variable}. Here we consider a DP method for generation of high-fidelity data. Development of accurate and efficient TRT solvers remains an active area of research, e.g. \cite{southworth2024one,southworth2025moment}. 
    
    The central problem of closure modeling is that dynamics at a microscale are subsumed into macroscale quantities using integral operators, which generically leads to nonlocal equations. Classically, closures have been derived using asymptotic arguments in limiting physical regimes. Recently data-driven methods for identifying closures have been proposed to extend the physical regimes of analytical closures, leveraging the availability of high-accuracy, high fidelity simulations. For linear neutral particle transport, entropy closure methods have been developed in \cite{hauck2011high,laiu2016positive,abdelmalik2023moment,schotthofer2025structure} and neural network-based closures are proposed in \cite{huang2022machine,porteous2023data,schotthofer2025structure}. A solver-in-the-loop framework for data-driven closure of gray TRT is presented in \cite{crilly2024learning}. Data-driven radiation diffusion-based closures have been developed in \cite{coale2019data,coale2024variable}. There has also been considerable progress made in data-driven closure identification for plasma physics \cite{huang2025machine,ingelsten2025data}, although a thorough review of this area is outside of the scope of this article, as our focus is on TRT.
    
    The Weak-form Sparse Identification of Nonlinear Dynamics algorithm (WSINDy) is a data-driven method that performs {\it equation discovery}, whereby governing equations are identified symbolically from data. The SINDy algorithm (on which WSINDy is based) and several variants are largely responsible for turning PDE identification from data into an accessible task, and accelerating research in this direction \cite{brunton2016discovering,rudy2017data,schaeffer2017learning}. SINDy proposes to discretize candidate operators and evaluate them at a given dataset (assumed to be an approximate solution to an underlying model), and then to fit the coefficients in front of each operator using sparse regression on the resulting linear system. Where classical methods of inference for differential equations, such as Markov-Chain Monte-Carlo (MCMC) and nonlinear least-squares, are typically {\it solver-in-the-loop} algorithms,
    SINDy does not perform forward solves of candidate models during the learning process\footnote{SINDy and variants are agnostic to how data is generated. Forward solvers are of course commonly used to {\it generate} data, as they are here, but this is distinct from the model learning process, which takes a {\it dataset} as input, and does not require a forward solver as input, in contrast to {\it solver-in-the-loop} algorithms.},  greatly accelerating model discovery, and allowing candidate equations to be enforced in a variety of ways. While the SINDy algorithm was originally posed in terms of the {\it strong form} of an equation, the key observation utilized in WSINDy is that a suitable {\it weak form} of the equation can be enforced instead. Originally introduced to reduce the effects of noise and non-smoothness  \cite{messenger2020weak,messenger2020weakpde}, WSINDy has been extended to many multiscale applications including mean-field equation discovery from particle data \cite{messenger2022learning}, identifying coarse-grained Hamiltonian systems \cite{messenger2024coarse}, and correcting analytical closures for equilibrating cold gases \cite{wang2025physics}. Weak-form equation learning is currently an active area of research  \cite{alves2022data,gurevich2024learning,russo2025streaming,TangLiaoKuskeEtAl2023JComputPhys,nicolaou2023data}; see the first author's SIAM News article \cite{messenger2024siam} and companion arXiv article (quoted therein) for references. 
    Despite avoiding forward solves in the learning process, it should be noted that SINDy-type algorithms do carry a computational bottleneck associated with sparse regression, which remains a combinatorially hard problem, although several recovery guarantees have recently been proven for strong and weak-form equation discovery \cite{he2022asymptotic,messenger2024asymptotic}.

    Of central importance in data-driven closures is well-posedness of resulting models, which often coincides with the enforcement of physical properties. It is in general an open problem to enforce well-posedness of data-driven models. Some preliminary work has been accomplished in \cite{kaptanoglu2021promoting}, which incorporates results from \cite{schlegel2015long} into SINDy to enforce stability in quadratic ODEs that obey an energy constraint. The problem of enforcing hyperbolicity in neural network closures was addressed in \cite{huang2023machine2,huang2023machine3} in the context of linear neutral particle transport, and in \cite{christlieb2025hyperbolic} for Boltzmann BGK. On the other hand, entropy closure methods are manifestly hyperbolic, and recent neural network-based entropy closures for radiation transport have been shown to be amenable to further structure preservation \cite{schotthofer2025structure}. 
 
    In this work, we develop a WSINDy-based framework to construct accurate moment closures to multifrequency TRT in optically thin media, a regime that does not admit a known analytical closure. To the best of the authors' knowledge, a data-driven method for identifying closures in optically thin multifrequency TRT, which avoids forward solves of candidate models in the learning process, has not been developed.
    Furthermore, we show that our learned closures can be parametrized over drive temperature and opacity level and extrapolated to produce adequate closure models in optically thin regimes, as characterized by a dimensionless Knudsen-like number. Outside of this region, in optically thick regimes, structurally different closures are expected, as well as an eventual transition to radiation diffusion. 
    We find closed systems in the form of hyperbolic balance laws for the radiation energy density $E$, radiation flux $F$, radiation-energy-weighted opacity $\sigma_E$, and temperature $T$, that implicitly provide a closure for the Eddington tensor $\CalE$ and model the nonlocal dependence on $\sigma$ through an auxiliary local equation. Under a simple change of variables, we find that many properties appear as convex constraints on model coefficients, including hyperbolicity, source stability, and preservation of equilibria, while other properties can be manifestly embedded into the model library, such as rotational symmetry. Lastly, our proposed WSINDy-based closures appear to be robust to numerical artifacts such as ray effects and particle noise, the former representing a significant corruption to the underlying physics in optically thin TRT simulations at commonly employed discretization levels (i.e., those used to benchmark the solver introduced in \cite{park2019multigroup}). From this, we conclude that our resulting closures may offer a unique alternative to mitigating numerical artifacts in future TRT solvers (see e.g.\ Fig.\ \ref{fig:high_low_res2}). More broadly, our results demonstrate that WSINDy is well-suited to augment the traditional search for kinetic closures by providing an interpretable data-driven framework in which physical properties can be enforced as convex constraints, and numerical artifacts do not corrupt learned closures.


In Section \ref{sec:TRT}, we review the TRT equations, including a description in Section \ref{sec:high_order_sims} of the TRT solver MuDDPaRT which we employ to generate high fidelity simulations, and related difficulties thereof. In Section \ref{sec:methods}, we describe our methodology, including the model class from which we will seek closures, physical constraints we impose, and the resulting WSINDy algorithm with constraints we use to identify closures. In Section \ref{sec:results}, we present results of our framework on a series of multifrequency 1D simulations with varying drive temperature and opacity. We focus on the ability of the closures to generate accurate predictions of moment quantities and to generalize over parameter space.

\section{Thermal radiation transport}\label{sec:TRT}
Consider the multifrequency time-dependent thermal radiation transport (TRT) and material temperature equations without scattering given by
\begin{subequations}\label{eq:ho}
\begin{align}\label{eq:ho-I}
  \frac{1}{c} \frac{\partial I}{\partial t} & = - \bOmega \cdot \nabla I - \sigma(\nu,T) I + \sigma(\nu,T)B(\nu,T), \\
  \rho c_V \frac{\partial T}{\partial t} & = \int_0^\infty \int_{S^d} \Big( \sigma(\nu,T) I  - \sigma(\nu,T)B(\nu,T) \Big)\d\bOmega\d\nu\label{eq:ho-T}.
\end{align}
\end{subequations}
Here, $I(\mathbf{x},\bOmega,\nu,t)$ is the specific radiation intensity for photons traveling in direction $\bOmega$ with frequency $\nu$, $T(\mathbf{x},t)$ is the material temperature, $c$ is the speed of light, $c_V$ is the specific heat (which in certain cases may depend on $T$), $\rho$ is the material density, and $\sigma(\nu,T)$ is the energy and temperature-dependent material opacity. The system \eqref{eq:ho} suffers from the curse of dimensionality, with $I$ living in a six-dimensional phase space plus time. A widely used 
framework for reducing the dimensionality is moment closure, whereby angle and frequency moments of $I$, together with temperature $T$, are modeled in place of the radiation intensity over all of phase space. Integrating over frequency and taking the zeroth and first angular moments, we get the moment system
\begin{subequations}\label{eq:lo}
\begin{align}
    \frac{\partial E}{\partial t} &= -\nabla \cdot \mathbf{F} - c\sigma_E(T)E + ac\sigma_P(T) T^4,\\
    \frac{1}{c}\frac{\partial \mathbf{F}}{\partial t} & = -\nabla\cdot (\mathcal{E} E) - \sigma_R(T) \mathbf F , \label{eq:moment-eqs-F}\\
    \rho c_V \frac{\partial T}{\partial t} & = -ac\sigma_P(T)T^4 + c\sigma_E(T) E,
\end{align}
\end{subequations}
where the radiation energy, flux, and Eddington tensor moments are defined by
\begin{equation}\label{eq:RT_moments}
  E \coloneqq \frac{1}{c} \iint I \d\nu \d\bOmega, \quad
  \mathbf{F}\coloneqq \iint  \mathbf{\bOmega} I \d\nu \d\bOmega, \quad    \mathcal{E}(I) \coloneqq c\frac{\iint \bOmega\bOmega I \d\nu\d\bOmega}{\iint I \d\nu\d\bOmega},
\end{equation}
and $\sigma_E,\sigma_P$, and $\sigma_R$ are the radiation-energy-, Planck-, and Rosseland-weighted opacities\footnote{Here $B$ is the black-body spectral radiance at temperature $T$ and frequency $\nu$,
\begin{equation}\label{eqB}
      B(\nu,T) = \frac{2h\nu^3}{c^2}\left( e^{h\nu/kT} - 1\right)^{-1},
\end{equation}
where $h$ and $k$ are the Planck and Boltzmann constants. This is derived from Planck's law that every physical body spontaneously and continuously emits electromagnetic radiation with spectral emissive power per unit area, per unit solid angle, per unit frequency, at temperature $T$, given by $B$. 
},
\begin{equation}\label{eq:collapse}
    \sigma_E(I,T) \coloneqq \frac{\iint \sigma I \d\nu\d\bOmega}{\iint I \d\nu\d\bOmega}, \hspace{4ex}
    \sigma_P(T) \coloneqq \frac{\int \sigma B \d\nu}{\int B \d\nu}, \hspace{4ex}
    \sigma_R(T) \coloneqq \frac{\int \frac{\partial B}{\partial T} \d\nu}{\int \frac{1}{\sigma}\frac{\partial B}{\partial T} \d\nu}.
\end{equation}
Conservation of total energy $e:=E+\rho c_VT$ is expressed by summing the evolution equations for $E$ and $T$.

For a given opacity $\sigma$, the general structure of \eqref{eq:lo} is a source-driven advection-reaction equation,  
\begin{equation}\label{eq:HBL_sec1}
    \partial_t \ubf = - \nabla \cdot \CalF(\ubf;I) + \CalS(\ubf;I), \qquad \ubf = (E,\Fbf,T),
\end{equation}
also known as a {\it Hyperbolic Balance Law}. Such equations are notoriously difficult to simulate and enforce stability thereof, as discussed below. Most importantly, \eqref{eq:HBL_sec1} is not closed due to dependence on $I$ through the Eddington tensor $\CalE$ and radiation-energy-weighted opacity $\sigma_E$. The predominant approach to solving this closure problem is to assume that the radiation is gray, such that $\sigma_E = \sigma_P$ and $\sigma_R$ are known functions of $T$ (or available in tabular form) and to argue that a diffusive approximation is valid to eliminate $\CalE$, whether via an isotropy assumption (no dependence on angle $\bOmega$) or by local thermal equilibrium. These approximations are reviewed in Appendix \ref{sec:known_closures}. These approaches can be extended to multifrequency radiation transport by discretizing in frequency space and then solving a coupled radiation diffusion system for all energy groups. 

In the optically {\it thin} regime, where the scattering media is rarefied and photon transport dominates over absorption and scattering, asymptotic assumptions leading to these well-known closures are not justified, and new methods must be developed to close \eqref{eq:lo}. This is the subject of the current paper. We aim to demonstrate that weak-form equation learning is especially well-suited to study this problem, given the available state-of-the-art methods for obtaining simulation data from the high fidelity system \eqref{eq:ho}. This hinges on the observation that WSINDy readily finds an auxiliary evolution equation for the radiation-energy-weighted opacity $\sigma_E E$, which closes the optically thin TRT dynamics in the form of a hyperbolic balance law. 

\subsection{High-order simulation method}\label{sec:high_order_sims}
We employ the moment-accelerated (high-order/low-order, HOLO \cite{chacon2017multiscale}) multi-group, deterministic particle code MuDDPaRT (Multi-Dimensional Deterministic Particle Ray-Tracer) to generate TRT data on \eqref{eq:ho} over one spatial dimension, referred to throughout as ``high fidelity data''. Introduced in \cite{park2019multigroup} for 1D and extended to multiple spatial dimensions in \cite{hammer2019multi}, MuDDPaRT discretizes \eqref{eq:ho} by representing the radiation intensity $I$ as a sum of Dirac-deltas over position-angle space, $I(\xbf,\bOmega,t) \approx \sum_{p=1}^Pw_p(t)\delta_{\vec{x}_p(t),\Omega_p(t)}$, where the particle equations of motions are given simply by $\dot{\vec{x}}_p = c\Omega_p$, $\dot{\Omega}_p = 0$, that is particles travel in straight lines with speed $c$. Particles are initialized at $k_i$ points in each spatial cell $C_i$, from each of which $M_\Omega$ particles are launched at angles $\{\hat{\Omega}_m\}_{m=1}^{M_\Omega}$. The weights $w_p$ are updated by integrating the characteristic equation for $I$ along each ray with implicit dependence on temperature at the next timestep, $T^{n+1}$. Each particle carries and updates all group intensities.

Each implicit time step is solved via a nonlinear high-order low-order (HOLO) iteration. On a high level, the low-order moment system \eqref{eq:lo} is utilized to obtain $T^{n+1}$, which is coupled to the high-order system through a consistency term $\qbf_c$ which represents the deviation from an isotropic radiation field,
\begin{equation} \label{eq:HOLO}
    \frac{1}{c} \frac{\partial \mathbf{F}}{\partial t} = -\frac{c}{3} \nabla E - \sigma_R(T) \mathbf{F} + c \qbf_c E.
\end{equation}
At each iteration, $\qbf_c$ is defined as the residual in terms of HO quantities (superscript $h$) at the previous iteration, 
\begin{equation}\label{eq:HOLO_q}
\qbf_c =  \left(\frac{1}{c} \frac{\partial \mathbf{F}^h}{\partial t} + \frac{c}{3} \nabla E^h + \sigma_R(T^h) \mathbf{F}^h\right)/ cE^h.
\end{equation}
The nonlinear LO system for $\{E,\mathbf{F},T\}$ is then solved according to the moment system, in particular to obtain an updated approximation of $T^{n+1}$. This new temperature and the particle trajectories are then used to update the radiation intensity $I$ equation, and a new iterate of $I$ is obtained by solving the HO equation for $I$. 
 This process of alternatingly solving HO and LO systems is iterated upon until the norm of the difference between low-order and high-order moments is below a user-specified tolerance, and implies discrete consistency between the high-order and low-order moment quantities to that tolerance before the next timestep.

\begin{figure}
\begin{tabular}{cc}
\hline\includegraphics[trim={0 0 20 0},clip,width=0.44\textwidth]{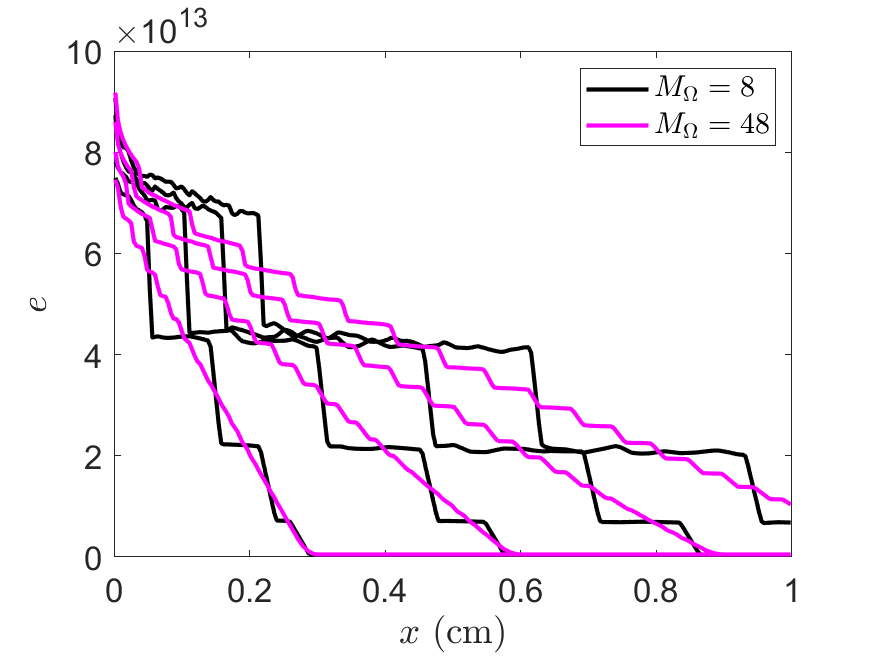} & 
\includegraphics[trim={0 0 20 0},clip,width=0.44\textwidth]{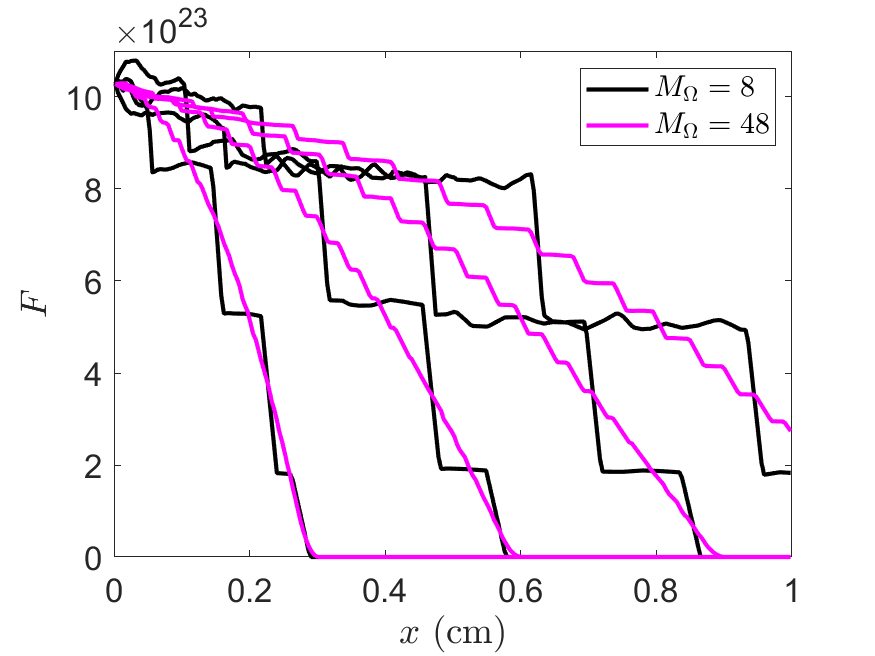} \\ 
\end{tabular}
\caption{Ray effects in MuDDPaRT simulations. Left and right plots show spatial profiles of the total energy $e$ and radiation flux $F$ at times $t\in \{1$e-$11s, 2$e-$11s, 3$e-$11s, 4$e-$11s\}$ for problem parameters $(\gamma,T_{in}^3)=(10^9,10^9)$, see \eqref{eq:paramgrid}. Black and magenta curves show simulations with $M_\Omega=8$ and $M_\Omega=48$ discrete polar angles, respectively. Increasing $M_\Omega$ serves to decrease jump sizes while increasing the total number of jumps.}
\label{fig:ray_effects}
\end{figure}

There are benefits and challenges associated with MuDDPaRT, as with any TRT solver, some of which motivate the current study. As a particle method, MuDDPaRT integrates the equations of motion more or less exactly in the radiation intensity equation, conserving energy locally and preventing diffusive numerical artifacts. MuDDPaRT is very efficient for multigroup calculations because each particle carries all energy groups, so tracking needs to only happen once to calculate all groups. MuDDPaRT is also asymptotic preserving by exactly capturing the diffusive limit. Compared with Implicit Monte-Carlo (IMC) \cite{fleck1971implicit}, deterministic particle tracking in MuDDPaRT avoids statistical noise, and is easily extensible to multigroup settings as the particle tracks can be reused for all frequency groups. The convergence rate of MuDDPaRT is empirically observed to be linear in the number of particles $\CalO(N_p^{-1})$, approaching quadratic in optically thin media, rather than $\CalO(N_p^{-1/2})$ as with Monte-Carlo methods, regardless of optical thickness \cite{park2019multigroup}.
On the other hand, unlike IMC particle methods, the discrete nature of angular discretizations in MuDDPaRT gives rise to {\it ray effects}, which manifest as spurious jumps in the low-order moment quantities. This is visualized in Figure \ref{fig:ray_effects}, where it can be seen that refining in angle space (from 8 to 48 polar angles) leads to visibly very different profiles for the total energy $e$ (defined in \eqref{eq:energycons}) and radiation flux $F$. Albeit less severe than statistical particle noise as in IMC simulations, deterministic particles also impart discretization errors that appear similar to statistical noise due to ray intersections with domain boundaries and cell boundaries (see e.g.\ $M_\Omega=8$ solution profiles in Fig.\ \ref{fig:ray_effects}).

\section{Methods}\label{sec:methods}
We seek a closure for the system of equations \eqref{eq:lo}, which necessitates finding functional expressions for the Eddington tensor $\CalE$ and radiation-energy-weighted opacity $\sigma_E(T)$ in terms of $E$, $\Fbf$, and $T$. We restrict the setting to one spatial dimension, as our focus is on finding closures for radiation transport with strong thermal dependence, anisotropy ($\CalE \neq \frac{c}{3}\Ibf$), and multifrequency opacities ($\partial_\nu \sigma \neq 0$). We also assume that $\sigma_P(T)$ and $\sigma_R(T)$ can be readily computed, since these are given analytically in terms of the opacity $\sigma$, which we assume to be known,
 and the Planck black-body spectral radiance $B$. Our primary considerations are then
\begin{enumerate}[label=(\Roman*)]
\item Preservation of known physics,
\item Stability and general well-posedness of the learned system of equations,
\item Feasibility of solving the resulting data-driven optimization problem,
\item Accurate standalone integration of $E$, $\Fbf$, and $T$ that well approximates moments and temperature as integrated from the true $I$.
\end{enumerate}
As we will see below, (I) tends to enable (II) and (III), which then leads to (IV). Below we describe the model class, optimization problem, and algorithms utilized to find closure models.

\subsection{Closure model class}


We found through preliminary applications of WSINDy that all diffusive and higher-order derivative terms were eliminated in the optically thin regime, reducing the search to a system of first-order PDEs. In particular, we found that $\sigma_E E$ approximately satisfies an evolution equation $\partial_t (\sigma_E E) = f(E,\Fbf,T,\sigma_EE)$ in weak form that is nonlinear and first-order in space and time derivatives, a description that is much more accurate and generalizes better than a direct functional relationship $\sigma_E = f(E,\Fbf,T)$, a gradient closure $\nabla \sigma_E = f(E,\Fbf,T)$, or a Laplacian closure $\Delta \sigma_E = f(E,\Fbf,T)$. Specifically, in our numerical experiments the direct functional relationship led to an inaccurate forward model, with further stability issues, while the gradient and Laplacian closures led to large weak-form equation residuals. This suggests that $\sigma_E$ depends nonlocally in space and time on other quantities, as expected from its definition \eqref{eq:collapse}. We choose to model $\CalE$ in terms of lower-order moments $E$ and $\Fbf$, as well as possibly coupling to temperature $T$, motivated in part by the iterative HOLO consistency method utilized in MuDDPaRT. As such, we utilize the evolution equation \eqref{eq:moment-eqs-F} for $\Fbf$ to project the divergence $\nabla \cdot (\CalE E)$ onto a library of terms that depend on $E$, $\Fbf$, and~$T$.  

We make an additional simplifying assumption that $\sigma_P(T) = \alpha / a cT^3$ for some $\alpha>0$, which holds in the examples below (as well as the Larsen and Marshak problems, which are widely used as benchmarks for forward solvers \cite{southworth2024one,park2019multigroup,hammer2019multi}), which is without loss of generality given our previous assumption that $\sigma_P(T)$ is available as a function of $T$ (implying that any $\sigma_P(T)$ can be inserted into the evolution equations for $E$ and $T$). Altogether, this leaves a system of equations of the form
\begin{subequations}\label{eq:lo_WS}
\begin{align}
    \partial_t E &= - \partial_x F - c\sigma_EE + \alpha T,\\
    \partial_t F &= - \partial_x p^F + q^F ,\\
    \partial_t T &= -\frac{\alpha}{\rho c_V} T + \frac{c}{\rho c_V}\sigma_E E ,\\
    \partial_t (\sigma_E E) &= - \partial_x p^\sigma  + q^\sigma ,
\end{align}
\end{subequations}
where we identify the closures for the fluxes $(p^F, p^\sigma)$ and source terms $(q^F, q^\sigma)$ through sparse constrained regression in weak form, as discussed below in Section \ref{sec:WSINDyoverview}. 

\subsection{Preserving physical properties}\label{sec:physical_properties}
We discuss several physical properties that any closure should obey. Each can be enforced naturally in the current sparse equation learning framework.

\subsubsection{Conservation of total energy}
The first physical property, which is preserved locally in MuDDPaRT, is conservation of total energy 
\[
    e := E +\rho c_V T.
\]
This follows from adding the evolution equations for $E$ and $T$, from which we get
\begin{equation}\label{eq:energycons}
\partial_t e +\nabla \cdot \Fbf = 0 .
\end{equation}
This is informative for closure identification, as the change of variables from $(E,\Fbf,T,\sigma_E E)$ to $(e,\Fbf,T,\sigma_E E)$ leads immediately to simplifications in the enforcement of hyperbolicity and source stability, as detailed below. From this we choose to identify closures that explicitly depend on $e$ rather than $E$, the latter obtainable from $e,\rho,c_V,T$. 

\subsubsection{Rotational symmetry}
The full TRT equations are rotationally symmetric, which can also be easily enforced in any closure at the library level. In one spatial dimension, this amounts to ensuring that reflections about the origin preserve the governing equations, or that the moment system is invariant to the transformation $(x,F)\to (-x,-F)$. This restricts the library significantly, in particular the closure system \eqref{eq:lo_WS} must obey
\begin{equation}\label{eq:refsym}
\begin{split}
p^F(\cdot,-F)&=p^F(\cdot,F), \quad q^F(\cdot,-F)=-q^F(\cdot,F) ,\\ p^\sigma(\cdot,-F) &= -p^\sigma(\cdot,F), \quad q^\sigma(\cdot,-F) = -q^\sigma(\cdot,F).
\end{split}
\end{equation}
This odd or even symmetry is enforced by restricting libraries to odd or even powers of $F$, respectively.

\subsubsection{Black-body equilibria}
A highly desirable physical quality of a closure model for TRT is that it preserves homogeneous equilibria of the kinetic system \eqref{eq:ho}. This is a relatively simple requirement that allows for numerical experiments to be conducted with initial and boundary conditions not seen by the data, as we demonstrate below. We can impose that known homogeneous equilibria of the kinetic system remain in equilibrium when passed into the closure model using equality constraints on the model coefficients. All homogeneous equilibrium TRT solutions are {\it black-body equilibria} and have the form $E = aT_*^4$,  $F=0$, and $\sigma_E E = a\sigma_p(T_*)T_*^4$ for $T_*$ constant. This can be enforced in the learned equations by selecting a grid of $N$ values $\Tbf = \{T_{\min},\dots,T_{\max}\}$, where for each value $T_i \in \Tbf$ we have an equilibrium $\ubf_i^* = (e,F,T,\sigma_EE)_i = (\rho c_V T_i + aT_i^4, 0, T_i, a\sigma_p(T_i)T_i^4)$. Suppose that for the $n$th equation we seek $J$ coefficients $\wbf^{(n)}\in \Rbb^J$ of terms from a library $\Lbb = \{\partial^{\alpha_1}f_1,\dots,\partial^{\alpha_J}f_J\}$, and that the indices $J_{source}^{(n)} \subset\{1,\dots,J\}$ correspond to terms without derivatives, $\partial^{\alpha_i} = \text{Id}$. We then build an equality constraint matrix $\Abf^{(n)} \in \Rbb^{N\times J}$ with entries 
\begin{equation}\label{eq:eqconstraints}
\Abf^{(n)}_{ij} = \begin{dcases} f_j(\ubf_i^*), & j\in J_{source}^{(n)} \\ 0, & \text{otherwise}\end{dcases},
\end{equation}
and impose the equality constraint $\Abf^{(n)}\wbf^{(n)} = 0$ in all least-squares solves. We find that this constraint is readily satisfied to high accuracy in the examples below.

Finally, it is known that black-body equilibria are stable to small perturbations\footnote{This can be seen by a maximum entropy argument, e.g. \cite{dubroca1999etude}}, hence we enforce linear stability of black-body equilibria with an analogous inequality constraint. Any numerical solver will impart small perturbations to the initial conditions, thus ensuring that an initially black-body equilibrium solution remains so requires that the model is linearly stable around such solutions. This is described in Section \ref{sec:stability}.

\subsection{Closure library}
We now discuss the library of terms $\Lbb$ we choose to represent our unknowns $(p^F, p^\sigma, q^F, q^\sigma)$. Despite the expected rational dependence of several terms on the state variables $E$ and $T$, experiments have demonstrated that polynomial terms are sufficient to model the dynamics in the examples below, and moreover that rational terms are removed from candidate models during sparsification (we retain several terms with inverse power laws in the equation for $\partial_t F$ to demonstrate this). The main challenge we face in designing a closure library is then the high dimensionality of the state space $(e,F,T,\sigma_E E)$. Rather than considering all monomial combinations of states, we generate function libraries using approximations from known relationships between variables as well as known physical properties of the TRT system. 

\subsubsection{Eddington closure}
We seek to close $\CalE E$ by directly modeling $\partial_x (\CalE E)$ as it appears in $\partial_t F$. First, we review several probabilistic ideas which align with previous closures as detailed in \cite{olson2000diffusion}. If we interpret $\bOmega$ as a random variable with distribution $\int I(\cdot,\nu)d\nu / \iint I(\bOmega,\nu)d \bOmega d\nu$, we note that 
\[\frac{1}{c}\CalE = \Ebb[\bOmega\bOmega^T] = \text{cov}[\bOmega]+\Ebb[\bOmega]\Ebb[\bOmega]^T = \text{cov}[\bOmega] + \left(\frac{\Fbf}{cE}\right)\left(\frac{\Fbf}{cE}\right)^T,\]
where $\Ebb[\cdot]$ and $\text{cov}[\cdot]$ denote the expectation and covariance operators. In the diffusive optically-thick limit, the radiation field is isotropic and $\bOmega$ is distributed uniformly over the sphere in $d$-dimensions leading to $\text{cov}[\bOmega] = \frac{1}{3}\Ibf_d$ and $\Fbf/cE \approx 0$. This is the well-known $P_1$ approximation \cite{olson2000diffusion}. In the opposite streaming limit, $\bOmega$ is $\delta$-distributed around $\frac{\Fbf}{cE}$ and hence $\text{cov}[\bOmega] = 0$. Another option is to model each component of $\bOmega$ as an independent and exponentially-distributed random variable $I(\Omega_i)\propto \exp(-\ell_i \Omega_i)$ over the interval $\Omega_i\in[-1,1]$ (with $\ell_i$ to be determined), in which case the relation
\[\Ebb[\Omega_i] = \frac{1}{\ell_i}-\coth(\ell_i) = \frac{F_i}{cE}\]
implies that the variance of each $\Omega_i$ satisfies
\[\text{var}[\Omega_i] = 1+\frac{1}{\ell_i^2}-\coth(\ell_i)^2 = 1+\frac{F_i}{cE}\left(\frac{2}{\ell_i}-\frac{F_i}{cE}\right).\]
Reflection symmetry (as discussed below) necessitates that $\CalE_{ii}$ $(= c\text{var}(\Omega_i)+\frac{F_i^2}{cE^2})$ is even with respect to $F_i$, which is given by eliminating odd powers of $F_i$ in the previous $\text{var}[\Omega_i]$ equation. One valid choice for $\ell_i$ is then $2/\ell_i = \alpha_i \frac{F_i}{cE}$ for constant $\alpha_i$, reducing the above to 
\[\text{var}[\Omega_i] = 1+(\alpha_i-1)\left(\frac{F_i}{cE}\right)^2\]
This aligns with the closures derived in \cite{kershaw1976flux,minerbo1978maximum}.

Other distributions can easily be considered to relate higher moments to lower moments. Using the distributions above, we expect $c\CalE E$ to depend on $(E,\Fbf)$ and possibly $T$, but not $\sigma_E E$, which informs our approximation of the flux as a general function of $(E,\Fbf)$. This becomes a general function of $(e,F,T)$ under the change of variables $E\to e$. For the 1D data considered below (with $u = \cos\phi$ being the projection of $\bOmega$ onto the $x$-axis), we thus seek an approximation to the angular variance
\[c\CalE E = c^2E \cdot\text{var}[u] + \frac{F^2}{E} \approx p^F (e,F,T),\]
where representation of $p^F$ as a polynomial suffices in the examples considered here. This is partially justified by considering the geometric series
\begin{equation}\label{eq:geom}
E^{-1} = E_0^{-1}\sum_{k=0}^\infty (1-E/E_0)^k,
\end{equation}
for $E_0$ satisfying $|E_0-E|<E_0$. We then model $\sigma_R(T)$ as a power law in $T$ given by $q^F(T)$, using our assumption that $\sigma_R(T)\propto T^{-3}$ (known to hold for opacity \eqref{eq:opacity} used in examples below, see Appendix \ref{sec:opacities}). The ansatz library $\Lbb^F$ is then defined by
\begin{equation}\label{eq:F_approx}
\begin{split}
\partial_t F &\approx -\partial_x p^F(e,F,T) + q^F(T) F \\&= \sum_{\substack{i+j+k=0\\ 0\leq i,j,k \leq p_{max}^F}}^{p_{tot}^F} \wbf^{p,F}_{ijk} \partial_x(e^iF^{2j}T^k) + \sum_{j=-p_{max}^F}^{p_{max}^F} \wbf^{q,F}_j T^j F =: \Lbb^F\wbf^F,
\end{split}
\end{equation}
for maximum total power $p_{tot}^F$ and maximum individual power $p_{max}^F$ defined by the user. Reflection symmetry is encoded in $p^F$ with even powers of $F$ and in $q^F$ by odd powers of $F$.

\subsubsection{Opacity closure}
Similarly, we parametrize the evolution equation for $\sigma_E E$ using
\begin{equation}\label{eq:sigma_approx}
\begin{split}
\partial_t (\sigma_E E) &= -\partial_x p^\sigma(F,T,\sigma_E E) + q^\sigma(e,T,\sigma_E E) \\
&= \sum_{\substack{i+j+k=0\\ 0\leq i,j,k\leq p_{max}^\sigma}}^{p_{tot}^\sigma} \wbf^{p,\sigma}_{ijk} \partial_x(F^{2i+1}T^j[\sigma_E E]^k) + \sum_{\substack{i+j+k=0\\ 0\leq i,j,k\leq p_{max}^\sigma}}^{p_{tot}^\sigma} \wbf^{q,\sigma}_{ijk} e^i T^j [\sigma_E E]^k =: \Lbb^\sigma\wbf^\sigma,
\end{split}
\end{equation}
with maximum total power $p_{tot}^\sigma$ and maximum individual power $p_{max}^\sigma$. The dependence of $p^\sigma$ on $(F,T,\sigma_E E)$ and $q^\sigma$ on $(e,T,\sigma_E E)$ is based on the following derivation of a leading-order model for $\sigma_E E$. Define the averaging operator $\lan f\ran := \iint f I d\bOmega d\nu/ \iint I d\bOmega d\nu = \iint f I d\bOmega d\nu/ cE$. Using that $\sigma=\sigma(\nu,T)$, we have 
\begin{align*}
    \partial_t (\sigma_E E) &= \frac{1}{c}\partial_t \iint \sigma(\nu,T(\xbf,t))I(\xbf,\bOmega,\nu,t)d\bOmega d\nu \\
    &= \frac{1}{c}\iint \Big[\partial_T\sigma(\nu,T(\xbf,t))\partial_t T I(\xbf,\bOmega,\nu,t) + \sigma(\nu,T(\xbf,t))\partial_t I(\xbf,b\Omega,\nu,t)\Big]d\bOmega d\nu \\
    &= \lan \partial_T \sigma\ran(\partial_t T) E + \iint \sigma(\nu,T) \Big(-\bOmega \cdot \nabla I - \sigma(\nu,T)(I-B)\Big)d\bOmega d\nu \\
    &=\frac{c\lan \partial_T \sigma\ran}{\rho c_V} \left(\sigma_E E - a\sigma_P(T)T^4\right) E   -c\nabla \cdot (E \lan \sigma \bOmega\ran )+ cE \nabla T\cdot \lan \partial_T \sigma \bOmega\ran - cE\lan \sigma\sigma\ran  +\iint \sigma^2 B d\bOmega d\nu,
\end{align*}
where we used $\sigma \nabla \cdot (\bOmega I) = \nabla \cdot (\sigma \bOmega I) - \partial_T \sigma \nabla T\cdot \bOmega I$ to arrive at the second and third terms on the right-hand side. Each average is then a single function or a product of functions. We then retain the leading order terms from the product averages, $\lan fg\ran\approx \lan f\ran \lan g\ran$, corresponding to an independence argument (in the probabilistic sense), or equivalently a low-rank approximation. Using that $\lan \bOmega\ran = F/cE$ and $\lan \sigma\ran = \sigma_E$, the equation then reduces to 
\begin{align*}
    \partial_t (\sigma_E E) &\approx \frac{c\lan \partial_T \sigma\ran}{\rho c_V} \left(\sigma_E E - a\sigma_P(T)T^4\right) E + \lan \partial_T \sigma \ran\nabla T\cdot \Fbf -\nabla \cdot ((E^{-1})\sigma_E E\Fbf) - cE^{-1}(\sigma_E E)^2 +q_0(T), 
\end{align*}
where we defined $q_0(T):=\iint \sigma^2 B d\bOmega d\nu$. For ease in identifying the equation in weak form, we use conservation of energy to transform $\nabla T\cdot \Fbf$ into the divergence of a flux plus source:
\[\nabla T\cdot \Fbf = \nabla \cdot(T\Fbf) -  T\nabla \cdot \Fbf = \nabla \cdot(T\Fbf) +  T\partial_t e \approx \nabla \cdot(T\Fbf) + T\partial_T e \partial_t T\]
\[ = \nabla \cdot(T\Fbf) + \frac{c T \partial_T e }{\rho c_V} \left(\sigma_E E - a\sigma_P(T)T^4\right), \]
where we used $\partial_te \approx \partial_T e \partial_tT$, assuming that temporal changes to energy approximately occur solely through temperature (holding e.g.\ in local thermal equilibrium, $e\approx aT^4 + \rho c_V T$). To find a leading-order polynomial representation of the dynamics, we invoke \eqref{eq:geom} again and retain the leading term $E_0^{-1}$. Together, this leads to the following simplified model in the form of a hyperbolic balance law:
\begin{align*}
    \partial_t (\sigma_E E) & \approx - \nabla \cdot \Big(E_0^{-1}\sigma_E E \Fbf -\lan \partial_T \sigma \ran T \Fbf\Big) + \\
    & + 
    \frac{c}{\rho c_V}\big(E\lan \partial_T \sigma\ran+ T \partial_T e\big) \big(\sigma_E E - a\sigma_P(T)T^4\big) - cE_0^{-1}(\sigma_E E)^2 +q_0(T) 
\end{align*}    
with functions $\lan \partial_T \sigma \ran$, $q_0$, and $\partial_T e$ to be determined. Assigning constants to each and using our assumption that $\sigma_P\propto T^{-3}$ gives us the leading-order polynomial representation in terms of $(e,\Fbf,T,\sigma_E E)$
\[\partial_t (\sigma_E E) \in \text{span}\Big\{\nabla \cdot(\sigma_E E\Fbf),\ \nabla \cdot(T\Fbf),\ e \sigma_E E,\ eT,\ T\sigma_E E,\ (\sigma_E E)^2,\ \tilde{q}_0(T) \Big\} , \]
where all purely temperature-dependent terms have been lumped into $\tilde{q}_0$. This justifies the term dependence in the expansion \eqref{eq:sigma_approx}. Since this holds to leading order, we directly enforce that the source terms $\{e \sigma_E E,\ eT,\ T\sigma_E E,\ (\sigma_E E)^2\}$ appear in the model (i.e., are not thresholded out during sparsification, see Appendix \ref{sec:WSINDyoverview}) with coefficients to be fit from data. We find that this greatly increases accuracy of the learned closure near the inflow wall in examples below.

\subsubsection{Final closed equations}
We will henceforth solve our closure identification problem in the context of the following {\it reflection symmetric total energy model}
\begin{subequations}\label{eq:ref_sym_energy_model}
\begin{align}
    \partial_t e  &= \wbf^0_1\partial_x F ,\\ 
    \partial_t F &= -\partial_x p^F(e,F,T;\wbf^{p,F}) + q^F(T;\wbf^{q,F})F ,\\
    \partial_t T &=  \wbf^0_2T + \wbf^0_3\sigma_E E ,\\
    \partial_t (\sigma_E E) &= -\partial_x p^\sigma (F,T,\sigma_E E; \wbf^{p,\sigma}) + q^\sigma(e,T,\sigma_E E; \wbf^{q,\sigma}),
\end{align}
\end{subequations}
with coefficients $\Wbf = (\wbf^{p,F},\wbf^{q,F},\wbf^{p,\sigma},\wbf^{q,\sigma}, \wbf^0)$ used to represent $(p^F,q^F,p^\sigma,q^\sigma)$ from the library $\Lbb = (\Lbb^F,\Lbb^\sigma)$ defined in \eqref{eq:F_approx} and \eqref{eq:sigma_approx}.  
The coefficients $\wbf^0$ are known analytically to equal $(-1,-\frac{\alpha}{\rho c_V},\frac{c}{\rho c_V})$, however we observe that WSINDy finds values that are better calibrated to the observed data due to the angular and energy discretizations utilized in MuDDPaRT. So we leave them as additional free parameters to be identified during the regression.

\subsection{Stability and Constraints}\label{sec:stability}

We now describe additions to our optimization problem that lead to a locally stable forward model. The closed system of equations \eqref{eq:ref_sym_energy_model} takes the form of a hyperbolic balance law, 
\begin{equation}\label{eq:HBL}
\partial_t \ubf + \partial_x \pbf(\ubf) = \qbf(\ubf),
\end{equation}
with $\ubf = (e,F,T,\sigma_E E)$. Well-posedness of such models has been studied in \cite{liu1979quasilinear,shizuta1985systems,chen1994hyperbolic,liu1999well,amadori2002global,hanouzet2003global,yong2004entropy}, with the hyperbolic case ($\qbf=0$) covered in \cite{friedrichs1971systems,godunov2025interesting}, and the textbook \cite{muller2013rational} offering a mathematical physics treatise under the umbrella of {\it extended thermodynamics}. 

\subsubsection{Local well-posedness}
Two necessary conditions for local well-posedness of models \eqref{eq:HBL} are the following:
\begin{enumerate}[label=(\roman*)]
\item (Hyperbolicity) $\partial_\ubf \pbf$ has real eigenvalues and is diagonalizable,
\item (Source stability) The ODE system $\frac{d\ubf}{dt} = \qbf(\ubf)$ has bounded trajectories.
\end{enumerate}
Using the method of characteristics, it is straightforward to see that the model gives rise to unbounded solutions if either (i) or (ii) is violated. However, (i) and (ii) are not sufficient for global well-posedness, which is often proven via entropy arguments. To the best of the authors' knowledge, global well-posedness for such models in general has only been proven for initial data close to equilibrium \cite{yong2004entropy}. 
Since global stability is in general an unsolved problem for practical cases (far from equilibrium), here we restrict our attention to ensuring that closed models are locally linearly stable by ensuring (i) weak hyperbolicity, that all eigenvalues of $\partial_\ubf \pbf$ are real, and a weakened version of (ii), that of linearly stable sources near equilibria, which are both necessary conditions for global stability as defined in \cite{yong2004entropy}. As discussed in the introduction, complete well-posedness in data-driven models is an open research topic.

\subsubsection{Hyperbolicity and linear source stability}
Hyperbolicity and linear source stability of the model \eqref{eq:ref_sym_energy_model} are each constraints on the eigenvalues of a state-dependent $n\times n$ Jacobian with $n=4$. This in general requires imposing nonlinear, nonconvex constraints at each grid point due to the $n$-degree polynomial dependence of the eigenvalues on the closure coefficients $\Wbf$. The total energy model \eqref{eq:ref_sym_energy_model}, together with the approximations made in \eqref{eq:F_approx} and \eqref{eq:sigma_approx}, enable both properties to be satisfied as convex constraints on $\Wbf$. Given the structure in \eqref{eq:ref_sym_energy_model}, we have 
\[\partial_\ubf \pbf = \begin{bmatrix} 0 & 1 & 0 & 0 \\ \partial_e p^F & \partial_F p^F & \partial_T p^F &0 \\ 0 & 0 & 0 & 0 \\ 0 & \partial_F p^\sigma & \partial_T p^\sigma & \partial_{\sigma_E E} p^\sigma  \end{bmatrix}, \]
which has eigenvalues $\lambda = \{0,\ d,\ \frac{1}{2}(b\pm \sqrt{b^2+4a})\}$ with $a=\partial_e p^F, b = \partial_F p^F, d=\partial_{\sigma_E E}p^\sigma$. A sufficient condition for $\lambda \in \Rbb$, and hence weak hyperbolicity, is then
\begin{equation}\label{eq:hyperbolicity}
\partial_e p^F\geq 0.
\end{equation}
This constraint is linear in the coefficients when $\partial_e p^F$ is evaluated at a grid of $(e,F,T)$ values. We thus follow the same protocol outlined above for enforcing black-body equilibria to enforce this constraint on a uniform grid in $(e,F,T)$.

Similarly, for linear source stability we need $\partial_\ubf \qbf$ to be negative semidefinite. We have 
\[\partial_\ubf \qbf = \begin{bmatrix} 0 & 0 & 0 & 0 \\ 0 & \partial_F q^F & \partial_T q^F &0 \\ 0 & 0 & -\frac{\alpha}{\rho c_V} & \frac{c}{\rho c_V} \\ \partial_e q^\sigma & 0 & \partial_T q^\sigma & \partial_{\sigma_E E} q^\sigma  \end{bmatrix}, \]
which has eigenvalues $\lambda = \{0, b, \frac{1}{2}[(h+d)\pm \sqrt{(h+d)^2-4(dh-fg)}]\}$ for $b=\partial_F q^F, h = \partial_{\sigma_E E} q^\sigma, d = -\alpha/\rho c_V, f = c/\rho c_V, g = \partial_T q^\sigma$. A sufficient set of linear constraints to ensure that eigenvalues of $\partial_\ubf \qbf$ have non-positive real part is the following 
\begin{subequations}\label{eq:sourcestability}
\begin{align}
\partial_F q^F&\leq 0 ,\\
\partial_T q^\sigma &\leq -\frac{\alpha}{c} \partial_{\sigma_E E}q^\sigma,\\ 
\partial_{\sigma_E E} q^\sigma &\leq \frac{\alpha}{\rho c_V}. 
\end{align}
\end{subequations}
We similarly enforce these on a grid of values. In this way, hyperbolicity and linear source stability follow from a small set of linear inequality constraints on the model coefficients $\Wbf$. In practice this is not enough to satisfy the constraints globally, thus we check {\it a posteriori} that each condition is satisfied over a given dataset. We find that sparsity promotion often leads to functional expressions that satisfy the given constraint over the parameter region considered. 

\subsection{WSINDy with constrained sequential thresholding}\label{sec:WSINDyoverview}

A brief review of WSINDy in the context of a single observed spatio-temporal field $u$ is included in Appendix \ref{app:WSINDy_overview}. Here we describe the extensions to WSINDy required to solve the constrained sparse regression problem for the free parameters $\Wbf = (\wbf^{p,F},\wbf^{q,F},\wbf^{p,\sigma},\wbf^{q,\sigma}, \wbf^0)$ in the reflection-symmetric energy model \eqref{eq:ref_sym_energy_model}. We represent the learned model as a sparse linear combination of terms from the library $\Lbb$ defined in \eqref{eq:F_approx} and \eqref{eq:sigma_approx}, which automatically embeds reflection symmetry. We use $p_{max}^F = p_{max}^\sigma=3$ and $p_{tot}^F=p_{tot}^\sigma=4$ to provide a balance between sparsity and accuracy, which leads to $\Lbb^F$ and $\Lbb^\sigma$ with 38 terms and 62 terms, respectively. This results in weak form matrices $\Gbf$ [see \eqref{eq:opt}] with conditions numbers $\kappa\approx \CalO(10^7)$. Larger libraries produced condition numbers too large to successfully find accurate sparse models, instead returning models with too few terms to accurately capture the dynamics.

Enforcing the model \eqref{eq:ref_sym_energy_model} in weak-form, subject to linear equality constraints to enforce black-body equilibria \eqref{eq:eqconstraints}, inequality constraints for linear stability of black-body equilibria \eqref{eq:sourcestability}, and inequality constraints for hyperbolicity \eqref{eq:hyperbolicity} leads to the optimization problem
\begin{equation}\label{eq:opt}
\min_\Wbf \|\Gbf(\Ubf)\Wbf - \Bbf(\Ubf)\|^2 +\lambda^2 \|\Wbf\|_0 \quad  \text{ s.t. } \quad \Abf\Wbf = \textbf{0},\ \Cbf\Wbf \leq \Dbf,
\end{equation} 
with simulation data $\Ubf = (\ebf,\Fbf,\Tbf,\Sbf)$ on the variables $\ubf = (e,F,T,\sigma_E E)$ and weak-form linear system $(\Gbf,\Bbf)$ given by the integration of terms in $\Lbb$ against test functions $\Psi = (\psi_1,\dots,\psi_K)$ [see Appendix \ref{app:WSINDy_overview} for details on the construction of ($\Gbf,\Bbf$)]. Here, $\Wbf$ is a $J\times n$ matrix of coefficients, $\Abf\in \Rbb^{M_E\times J}$ is a matrix representing equality constraints \eqref{eq:eqconstraints}, and $\Cbf\in \Rbb^{M_I\times J}, \Dbf\in \Rbb^{M_I\times n}$ are matrices representing inequality constraints \eqref{eq:hyperbolicity}, \eqref{eq:sourcestability}. We solve \eqref{eq:opt} using sequentially thresholded quadratic programming. For the thresholding part of the algorithm, we use the MSTLS method introduced in \cite{messenger2020weakpde} (and reviewed in Appendix \ref{app:WSINDy_overview}), with inner least-squares solves replaced with quadratic programs to enforce the convex constraints. We use an interior-point algorithm with 1000 maximum iterations and constraint and optimality tolerances of $1$e-$10$. In our case, $\Dbf\geq \textbf{0}$ entry-wise, so $\Wbf=0$ is always a feasible point. Feasibility of nontrivial $\Wbf$ requires $\Abf$ to have incomplete column rank, hence we enforce black-body equilibria on a small grid of temperatures $\Tbf = \texttt{linspace}(0,4T_{max},7)$, where $T_{max}$ is the maximum temperature value observed in the high-fidelity data. Linear source stability is enforced at $\Tbf$ together with 20 solution values along the inflow boundary, $\ubf(0,t_i)$ for $i=1,\dots,20$, where the data exhibits the greatest particle errors. We evaluate the hyperbolicity constraint $\partial_e p^F\geq 0$ on a grid $\texttt{linspace}(0,e_{max},p_{max}^F)\times \texttt{linspace}(0,F_{max},p_{max}^F+1) \times \texttt{linspace}(0,T_{max},p_{max}^F+1)$ where $(e_{max},F_{max})$ are the maximum values observed for $(e,|F|)$ and $p_{max}^F=3$ throughout. This leads to $47$ total equality constraints and $149$ inequality constraints. We find this to be sufficient but note that we have avoided including an excessive number of constraints in order to expedite optimization, which was performed in MATLAB on a laptop with 32 Gb of RAM. For larger problems with more computational resources, more constraints can easily be added.  

To enforce equations in the weak form, we use rectified piecewise polynomial test functions $\varphi(v) = (1-(v/a)^2)_+^p$ separately along the space and time dimensions (where $[x]_+ := \max\{x,0\}$), implemented using the convolutional approach review in Appendix \ref{app:WSINDy_overview}. Such test functions were shown in \cite{messenger2020weak} to offer high-accuracy coefficient recovery from noiseless data for sufficiently large $p$. An algorithm to select hyperparameters $(p,a)$ from noisy data was introduced in \cite{messenger2020weakpde}. We use this method with hyperparameters $(\tau,\hat{\tau}) = (10^{-4},6)$, where $\tau>0$ and $\hat{\tau}>0$ set the decay rates of the test function in real and Fourier space, respectively, and subsequently define $(p,a)$ (see the discussion following \eqref{septest} for details). Thus from $(\tau,\hat{\tau})$ we define a total of 16 test function parameters, $(p_x,p_t)$ and $(a_x,a_t)$ for each of the four equations \eqref{eq:ref_sym_energy_model}, for each dataset. We find that $p \in \{3,4,5\}$, $a_x\approx L_x/15$, and $a_t\approx L_t/4$, where $L_i$ is the length of the domain along the given axis. We note that this may be counterintuitive from the perspective of finite elements, where typically $p\in \{1,2\}$ and $a$ is on the order of the meshwidth. In WSINDy, larger support $a$ helps to remove the effects of noise and other corruptions, while larger $p$ leads to more accurate numerical integration over smooth regions of the data. In addition, where previous work utilized even larger powers $p$ to increase coefficient accuracy \cite{messenger2020weak,messenger2020weakpde}, we did not see significant gains from this due to discontinuities in the data from ray effects limiting accuracy across the entire spacetime domain. 

\subsection{Parametrized closures}

To promote generalization, we identify closures simultaneously from a small set of high fidelity simulations $\Ubb = \{\Ubf^{(1)},\dots,\Ubf^{(P)}\}$, each of which is simulated from a different value $\pmb{\mu}^{(p)}$ of the physical parameters involved [in the cases below, $\pmb{\mu} = (T_{in},\gamma)$]. To enforce the same closure across each simulation, we use group sparsity in combination with MSTLS as outlined above. In this case we replace thresholding constraints of the form 
\[\ell^{(p)}_i(\lambda)\leq |\wbf^{(p)}_i|\leq u^{(p)}_i(\lambda),\]
[see \eqref{MSTLSbnds}] with their analogous $\ell_1$-norm concatenations,
\[\sum_{p=1}^P\ell_i^{(p)}(\lambda)\leq \|\wbf^{(p)}_i\|_1 \leq \sum_{p=1}^Pu_i^{(p)}(\lambda).\]
In this way, the quadratic programming solves for each simulation can occur in parallel, and the coefficients are coupled together only through the thresholding constraints. The MSTLS loss \eqref{lossfcn} is then taken to be the sum of the respective losses over each simulation, 
\[\CalL^{group}(\lambda;\Ubb):= \sum_{p=1}^P\CalL(\lambda;\Ubf^{(p)}).\]
Once a set of $P$ coefficients $\Wbf^{(p)}$, $p=1,\dots,P$, is found for each simulation, a suitable interpolation scheme must be utilized to interpolate and extrapolate the closure at new parameter values $\pmb{\mu}$. While interpolation schemes in parametrized PDE discovery and reduced-order modeling is an active area of research \cite{fries2022lasdi,nicolaou2023data,tran2024weak}, in the cases below we find that coefficients can be extrapolated log-linearly, providing another layer of interpretability to the closure. 

\section{Results}\label{sec:results}
We aim to demonstrate that the proposed constrained WSINDy approach is an advantageous framework for identifying closures in optically thin TRT and that the resulting hyperbolic balance laws \eqref{eq:ref_sym_energy_model} generalize well over a range of problem parameters. We focus on the case of learning a closure from a small set of high-fidelity simulations and then parametrizing the closure coefficients so that the closure may be tested on new problem parameters. We note that the proposed framework successfully solves the simpler problems of (a) learning a closure model from a single simulation that can reconstruct the given data (indeed, Fig.\ \ref{fig:highres_closure} shows this, as well as the closure performance on the training data, see Fig.\ \ref{fig:extrap}), and (b) interpolating the learned closure model to produce solutions at parameters inside the convex hull of the training points, but not equal to any of the training points (verified, not shown). Hence, for conciseness we omit these studies and present only results focused on extrapolation. By studying the extrapolation capabilities of the parametrized closure, we identify regions of parameter space where the closure is expected to perform well, and regions where transition to a diffusive closure is expected (an investigation we reserve for future work). Lastly, we quantify the performance of the closure in light of ray effects in the data, revealing that closure solutions are devoid of ray effects, suggesting that learned closures, as we have presented here, may provide a way to mitigate ray effects in future solvers.

We restrict our examples to a representative one-dimensional multifrequency problem consisting of radiation from a ``hot'' material flowing into an initially ``cold'' material, both initially black bodies in thermal equilibrium.\footnote{In what follows we report temperature values in units of electron volts (eV), which implies the Boltzmann constant is unity, $k=1$, lengths in centimeters (cm), and time in seconds (s).} The temperatures of the hot and cold materials are $T_{in}\sim 10^3\text{eV}$ and $T_o\sim 1\text{eV}$, the former providing boundary conditions while the latter provides initial conditions. We consider the opacity from the classical Larsen problem \cite{larsen1988grey},
\begin{equation}\label{eq:opacity}
\sigma(\nu,T) = \frac{\gamma}{(h\nu)^3}(1-e^{-h\nu/k T}),
\end{equation}
which leads to $\alpha = \frac{15 ac}{\pi^4 k^3}\gamma = \frac{60\sigma_{SB}}{\pi^4k^3}\gamma$ in \eqref{eq:lo_WS}. We examine learned closures for variable inputs $(\gamma, T_{in})$ and fixed $(T_o,\rho,c_V)$, using high-fidelity simulations of the TRT dynamics for $x\in [0,4\text{cm}]$ over 1024 uniformly-spaced gridpoints and inflow boundaries at $x=0$. Data on $(e,F,T,\sigma_E E)$ is collected for 200 uniform timesteps at $\Delta t = 1$e-$12$. In order to demonstrate that the proposed approach performs well in the presence of ray effects, we predominantly focus on high-fidelity data obtained from MuDDPaRT at the low-resolution angular discretization $M_\Omega=8$. A comparison to $M_\Omega = 48$ data is included in Section \ref{sec:res_effects}.

We show that the resulting learned closure generalizes to new spatiotemporal domains and parameter values, a central benefit of the interpretability of equation learning compared to black-box machine learning methods.
We exhibit this by training on the simulation domain $(x,t) \in [0,2]\times [0,1$e-$10]$ and simulating the learned closure models over $(x,t) \in [0,4]\times [0,2$e-$10]$. Figure \ref{fig:extrap_33} visualizes this for the problem parameters $(\gamma,T_{in}^3) = (10^9,10^9)$, where the training domain overlays the simulation data in the top row. To exhibit reflection symmetry of learned models, we simulate closure models under the transformation $(x,F)\to (-x,-F)$ on the initial and boundary conditions, and compare to the high-fidelity data under the same transformation. For forward simulation of closures, we use a 5th-order WENO discretization in space \cite{jiang1996efficient} with twice the meshwidth of the high fidelity data and the adaptive Runge-Kutta 4/5 method in time. We use Dirichlet boundaries for $(e,F,\sigma_E E)$ at the inflow boundary given by the high-fidelity data fit to a quadratic function in time\footnote{In general, boundary conditions represent an additional closure problem. We reserve this for future work.}. Two ghost nodes are used to extend the WENO stencil at each boundary by the solution value at the boundary.

We use the following three relative error metrics to quantify the accuracy of learned closures:
\begin{equation}\label{eq:metrics}
    \text{err}_{L_1}(u) := \frac{\sum_{i=1,j=1}^{N_x,N_t}| \widehat{\Ubf}_{ij}-\Ubf_{ij}|}{\sum_{i=1,j=1}^{N_x,N_t}|\Ubf_{ij}|}, \quad     
    \text{err}_{{L_1},j}(u) := \frac{\sum_{i=1}^{N_x}| \widehat{\Ubf}_{ij}-\Ubf_{ij}|}{\sum_{i=1}^{N_x}|\Ubf_{ij}|}, \quad 
    \text{err}_{Int,j}(u) := \frac{\left\vert\sum_{i=1}^{N_x}(\widehat{\Ubf}_{ij}-\Ubf_{ij})\right\vert}{\left\vert\sum_{i=1}^{N_x}\Ubf_{ij}\right\vert},
\end{equation}
where $\widehat{\Ubf}$ and $\Ubf$ represent the simulated closure solution and high-fidelity data, respectively, for approximating the state variable $ u\in \{e,F,T,\sigma_E E\}$. The grid resolution used in simulations of the closure model is $N_x = 512$ points in space and $N_t = 200$ points in time (the high-fidelity data is interpolated onto this grid, which simply amounts to subsampling in space). The metric $\text{err}_{L_1}$ is the relative $L_1$ norm over spacetime, while $\text{err}_{L_1,j}$ is the relative $L_1$ norm over space at time $t=t_j=j\Delta t$. As discussed below, these metrics can be misleading due to ray effects, so we also include the relative {\it integrated error} $\text{err}_{Int,j}$ as a weaker metric, which due to sign definiteness of our examples is equivalent to $\left\vert\|\widehat{\Ubf}_{:,j}\|_1-\|\Ubf_{:,j}\|_1\right\vert / \|\Ubf_{:,j}\|_1$, the relative difference between the spatial $L_1$ norms of the closure solution and high-fidelity data, at time $t=j\Delta t$. These and other weak metrics are used in weak convergence studies of numerical solutions to stochastic differential equations (SDEs) \cite{higham2001algorithmic}, for example.

\subsection{Problem parameters and associated dimensionless number $\kappa_L$}

To characterize the performance of our closures, we define a Knudsen number-type dimensionless number $\kappa_L$ as the ratio of $1/\sigma_P(T_o;T_{in})$, a measure of the mean-free path, to a characteristic length $L$, where  
\begin{equation}
\sigma_P(T_o;T_{in}) :=   \frac{\int_0^\infty \sigma(\nu,T_o)B(\nu,T_{in})d\nu}{\int_0^\infty B(\nu,T_{in})d \nu}.
\end{equation}
The quantity $\sigma_P(T_o;T_{in})$ is the opacity of the ambient material at $T_o$ averaged over the density of incoming photons near the wall, in other words the inverse of the average optical depth near the wall. We then define $\kappa_L$ by 
\begin{equation}\label{eq:kappaL}
\kappa_L := (L\sigma_P(T_o;T_{in}))^{-1},
\end{equation}
where we use the domain size $L=4\text{cm}$ as our characteristic length. For opacities of the form $\sigma(\nu,T) = \frac{\gamma}{(h\nu)^3} f\left(\frac{h \nu}{k T}\right)$ for some function $f:[0,\infty)\to[0,\infty)$, we have that 
\[\sigma_P(T_o;T_{in}) = \frac{15 \gamma }{\pi^4 (kT_{in})^3} \int_0^\infty \frac{f\left(\frac{T_{in}}{T_o}x\right)}{\exp(x)-1}dx, \]
which for the specific case \eqref{eq:opacity} leads to, when $T_{in}\gg T_o$, 
\begin{equation}\label{eq:sigma_PTio_approx} 
\sigma_P(T_o;T_{in})\approx \frac{15\gamma}{\pi^4 (kT_{in})^3} \log\left(\frac{T_{in}}{T_o}\right)
\end{equation}
(see Appendix \ref{sec:opacities}). The problem becomes optically thick when $\kappa_L\ll 1$, or when the number of mean-free paths per $L$ is very large.



Given the approximate power laws $\kappa_L  \sim T_{in}^3/ \gamma$ derived from \eqref{eq:sigma_PTio_approx}, we examine closures as $\gamma$ and $T_{in}^3$ are varied over two orders of magnitude centered at the nominal values from the optically thin region of the Larsen problem. That is, we consider 
\begin{equation}\label{eq:paramgrid}
    (\gamma,T_{in}^3)\in \{10^8,10^{8.5},10^9,10^{9.5},10^{10}\}^2.
\end{equation}
We use data from a subset of parameters to learn the functional form of the closure, and then interpolate the coefficient values using bivariate splines in logarithmic scale to sample the closure at new parameter points, having observed {\it a posteriori} that the coefficients exhibit an approximately log-linear parameter dependence.

\subsection{An optically thin closure}

For our optically thin closure, we train on simulations at the nine most optically thin parameter values (equivalently the largest values of $\kappa_L$)  $(\gamma,T_{in}^3) \in \{10^8,10^{8.5},10^9\}\times\{10^9,10^{9.5},10^{10}\}$, for the purpose of reproducing the training dynamics as well as assessing the extent to which the closure can be used on the remainder of the parameter grid. The resulting closures for $F$ and $\sigma_E E$ have the form
\begin{equation}\label{eq:extrap_closures} 
\begin{dcases} 
\partial_t F &= \partial_x\left( \wbf^{p,F}_1  e+ \wbf^{p,F}_2eF^2 + \wbf^{p,F}_3F^4 \right),\\ 
\partial_t (\sigma_E E) &=  \partial_x\left(\wbf^{p,\sigma}_1  F +\wbf^{p,\sigma}_2 (F\sigma_E E) +\wbf^{p,\sigma}_3 (F (\sigma_E E)^3) + \wbf^{p,\sigma}_4 F^3 +\wbf^{p,\sigma}_5  F^3 \sigma_E E\right) \\ 
&\qquad +\wbf^{q,\sigma}_1 (\sigma_E E)^2 +\wbf^{q,\sigma}_2  T \sigma_E E +\wbf^{q,\sigma}_3 e\sigma_E E + \wbf^{q,\sigma}_4 eT .
\end{dcases}
\end{equation}
The Eddington closure implied from $\partial_t F$ takes a different form from those typically found in the literature due to the implicit temperature dependence incorporated in the total energy $e$. Dependence on $F^4$ indicates that the underlying radiation intensity $I$  interpreted as a probability over angle space lies outside of the commonly used exponential family, which includes the delta, uniform, and exponential distributions discussed above. We also find that the sink term $\sigma_R(T) F$ is lower-order in the optically thin regime and hence not identified, possibly due again to the implicit temperature dependence already captured in the flux. For the opacity closure, the dynamics in the bulk are predominantly given by the flux $p^\sigma \sim F$, while temperature-driven source terms effectively capture sharp increases in opacity near the wall. 
\begin{figure}
\begin{tabular}{@{}c|@{}c@{}}
    \includegraphics[trim={0 50 20 25},clip,width=0.25\textwidth]{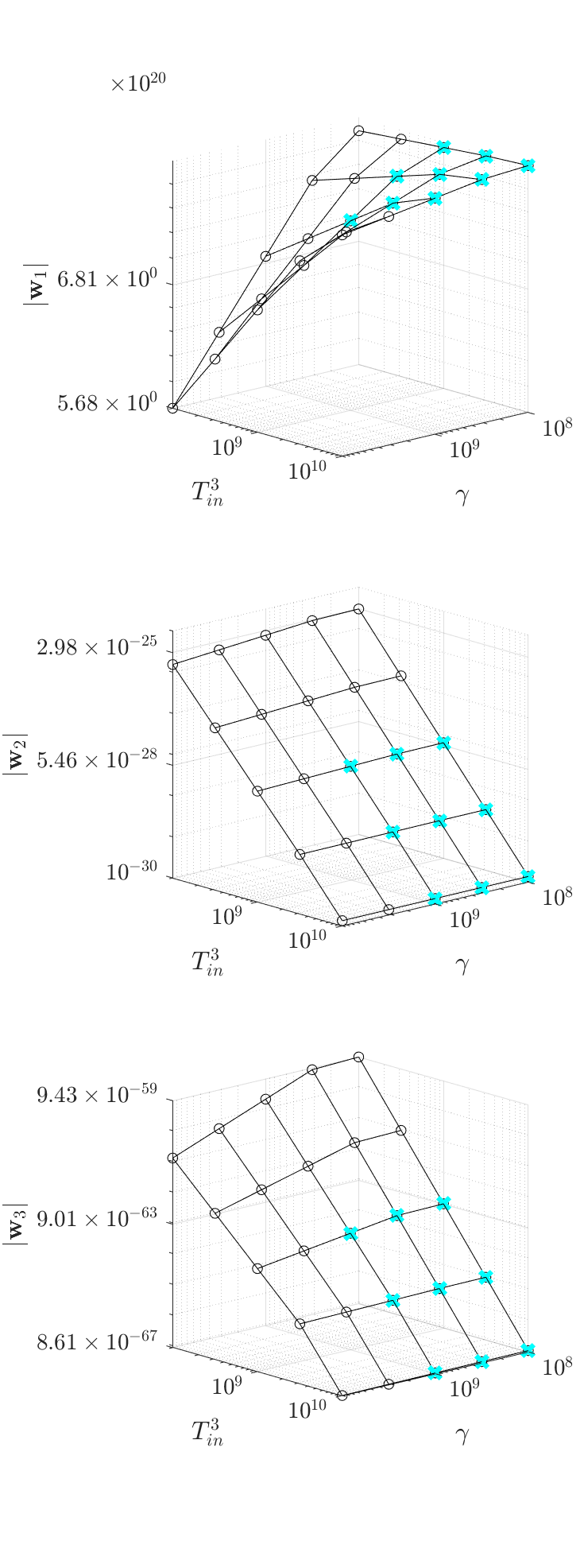} & 
    \includegraphics[trim={70 75 90 80},clip,width=0.75\textwidth]{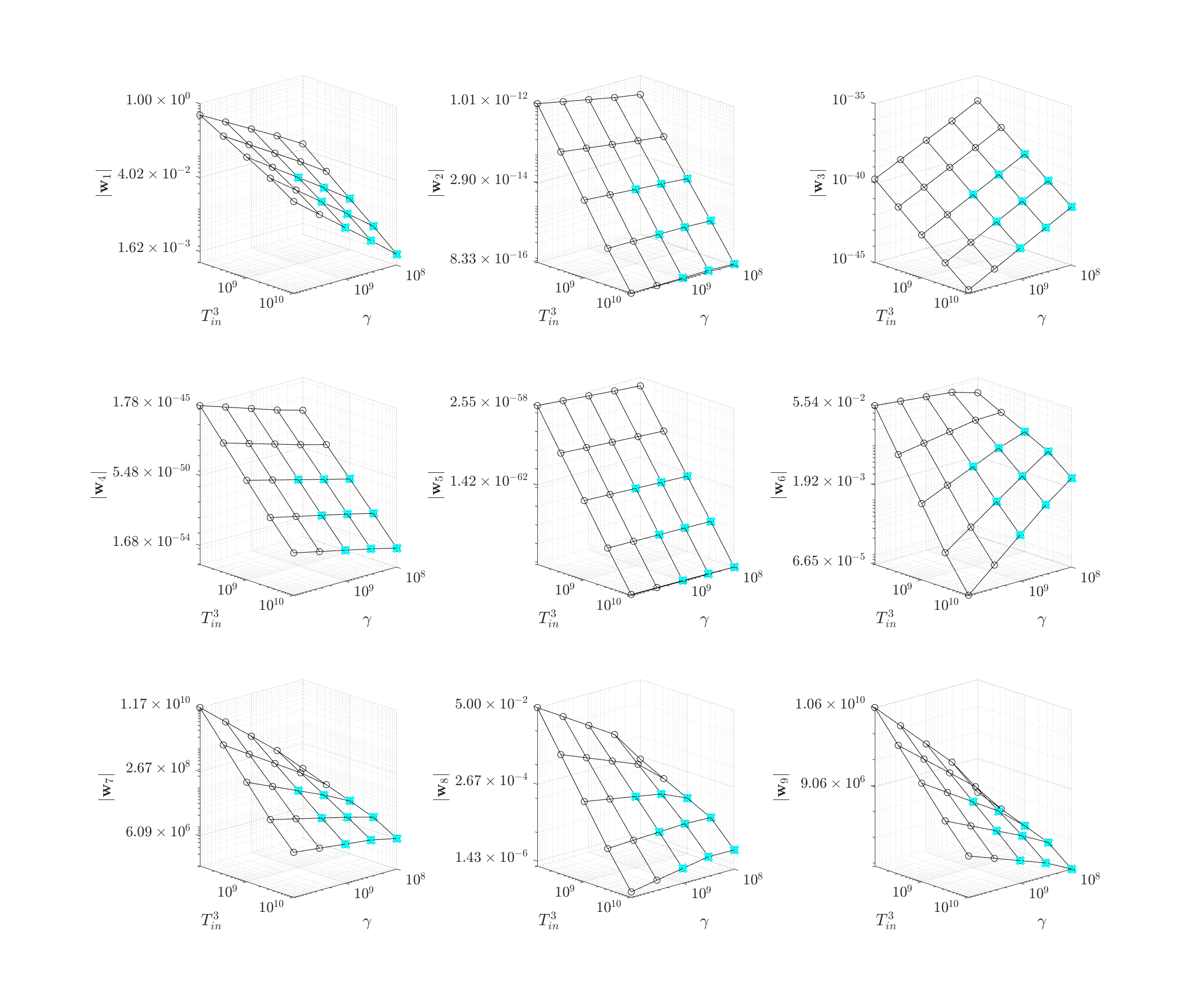}
\end{tabular}
\caption{Coefficient values for closure model \eqref{eq:extrap_closures} plotted in log space. Left and right show $\partial_t F$ and $\partial_t (\sigma_E E)$ coefficients. Cyan dots indicate coefficient values found for training simulations. The rest of the values are extrapolated from the training values in a log-linear fashion.}
\label{fig:coeff_grid}
\end{figure}

Coefficients $\Wbf$ learned for the nine training simulations are visualized in Figure \ref{fig:coeff_grid}, along with their extrapolation to the remaining parameter points. Using the fact that coefficients have the same sign across training simulations (this can easily be enforced, but was not in the closures presented), we plot the absolute values of coefficients logarithmically across the $(\gamma,T_{in}^3)$ plane. From this it can readily be seen that most coefficients lie approximately on a plane, justifying log-linear extrapolation. Hence, we model each set of coefficients $\wbf_i$ using the parametric relation
\begin{equation}\label{eq:loglinrelation}
\wbf_i \approx \wbf_{i,0} \gamma^{\eta^\gamma_i} (T_{in}^3)^{\eta^T_i},
\end{equation}
which is use to extrapolate the closure in the next Section. The reference coefficient $\wbf_{i,0}$ and powers $(\eta^\gamma_i,\eta^T_i)$ are given in Table \ref{tab:coeffs}. From this we see that small changes in $T_{in}$ lead to drastic changes in several of the coefficients, while changes in $\gamma$ have less severe impact. The coefficient $\wbf_1^{p,F}$ of $\partial_x e$ stays relatively constant, changing by approximately $10\%$ as $\gamma$ is varied over two orders of magnitude. This is not suspect, as its base coefficient approximates the correct isotropic limit $-c^2$ to better than $95\%$ accuracy. In contrast, the coefficient $\wbf_5^{p,\sigma}$ changes by 8 orders of magnitude as $T_{in}^3$ is varied over the same range. The relative errors of the log-linear coefficient fit (last column) indicate that the trend very accurately matches the coefficient data. 

\begin{table}
\centering
    \begin{tabular}{c|c|ccc|c}
         & Term&  $\wbf_{i,0}$ & $\eta^T$ & $\eta^\gamma$ & $r^2$ \\ \hline
        $\wbf_1^{p,F}$ & $\partial_x e$& -8.17e+20 & 1.56e-02 & -1.87e-02 & 2.31e-04\\
        $\wbf_2^{p,F}$ & $\partial_x(eF^2)$ &6.80e-04 & -2.65 & -2.02e-02 & 2.86e-04 \\
        $\wbf_3^{p,F}$ & $\partial_x(F^4)$ & -9.65e-27 & -3.91 & -1.01e-01 & 8.11e-04 \\ \hline 
        $\wbf_1^{p,\sigma}$ &$\partial_xF$ &-3.49e-02 & -7.06e-01 & 7.19e-01 & 8.51e-03\\
        $\wbf_2^{p,\sigma}$ &$\partial_x(F\sigma_E E)$ & -4.01e-02 & -1.39 & 2.62e-02 & 8.66e-04\\
        $\wbf_3^{p,\sigma}$ &$\partial_x(F(\sigma_E E)^3)$ &8.95e-05 & -2.38 & -1.69 & 3.30e-04\\
        $\wbf_4^{p,\sigma}$ &$\partial_x(F^3)$ & 7.24e-26 & -3.51 & 8.03e-01 & 3.33e-04  \\
        $\wbf_5^{p,\sigma}$ &$\partial_x(F^3\sigma_E E)$ &1.70e-26 & -4.08 & 5.72e-02 & 2.11e-04 \\
        $\wbf_1^{q,\sigma}$ &$(\sigma_E E)^2$ & -8.60e+09 & -7.67e-01 & -6.10e-01 & 2.76e-02\\
        $\wbf_2^{q,\sigma}$& $T\sigma_E E$  &1.82e+11 & -7.67e-01 & 3.90e-01 & 9.69e-03 \\
        $\wbf_3^{q,\sigma}$ & $e\sigma_E E$ &3.07e+07 & -1.38 & 1.07e-01 & 2.45e-02 \\
        $\wbf_4^{q,\sigma}$ & $eT$ & -6.50e+08 & -1.38 & 1.11 & 2.29e-02
       \end{tabular}
       \caption{Log-linear extrapolation parameters for each coefficients in the closure model \eqref{eq:extrap_closures}. The last column shows the relative $\ell_2$ error of the coefficient fit.}
   \label{tab:coeffs}
\end{table}

\subsection{Extrapolating in parameter space}\label{sec:extrap}

\begin{figure}
\begin{tabular}{@{}c@{}|c@{}}
\begin{minipage}{0.6\textwidth}

\centering
        \begin{subfigure}[b]{0.49\textwidth}
            \centering
                \includegraphics[trim={0 0 0 0},clip,width=1\textwidth]{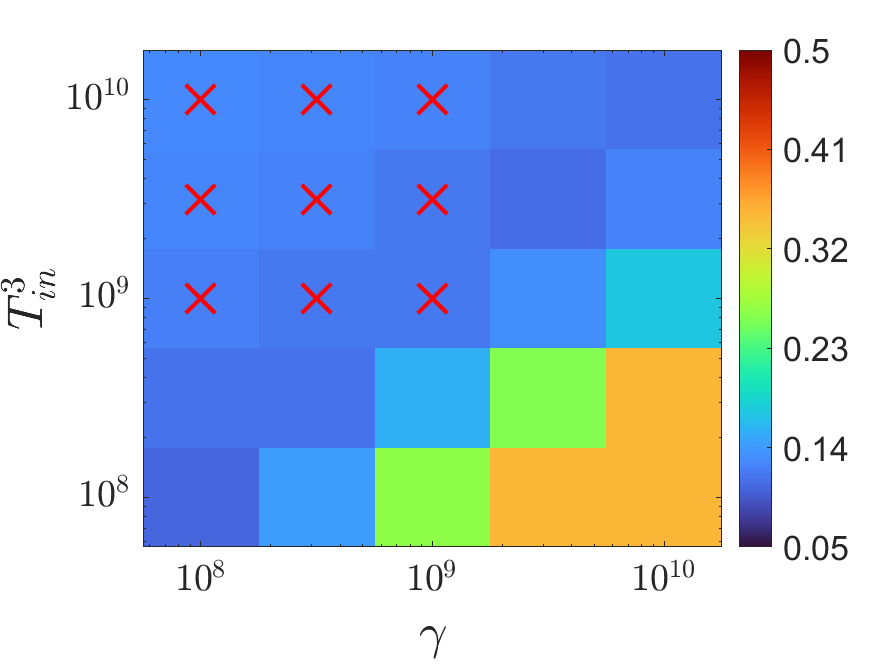}
            \caption[Network2]%
            {{\small $\|e - \widehat{e}\|_1/\|e\|_1$}}    
            \label{fig:mean and std of net14}
        \end{subfigure}
        \begin{subfigure}[b]{0.49\textwidth}  
            \centering 
                \includegraphics[trim={0 0 0 0},clip,width=1\textwidth]{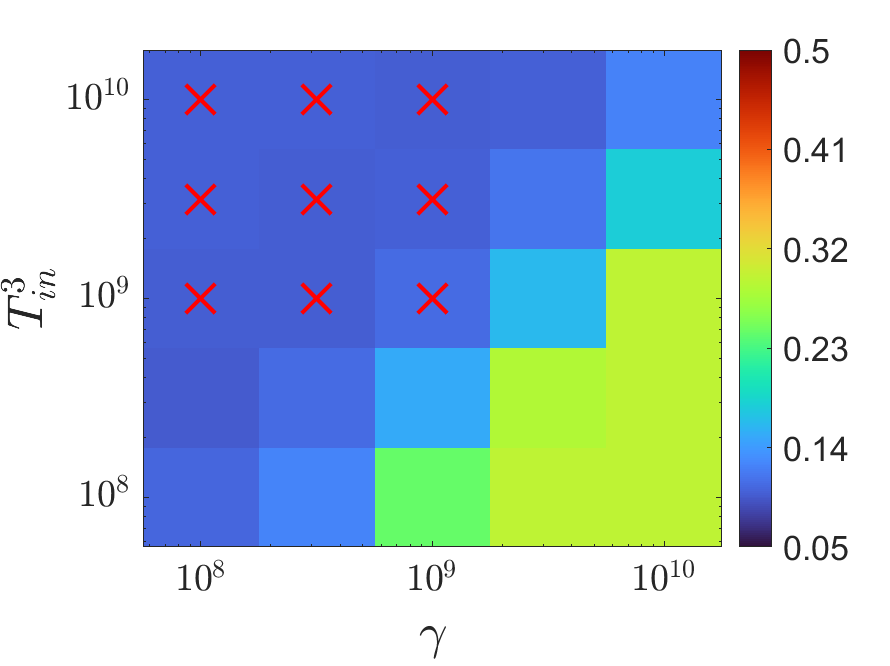}
            \caption[]%
            {{\small $\|F - \widehat{F}\|_1/\|F\|_1$}}    
            \label{fig:mean and std of net24}
        \end{subfigure}
        \vskip\baselineskip
        \begin{subfigure}[b]{0.49\textwidth}   
            \centering 
                \includegraphics[trim={0 0 0 0},clip,width=1\textwidth]{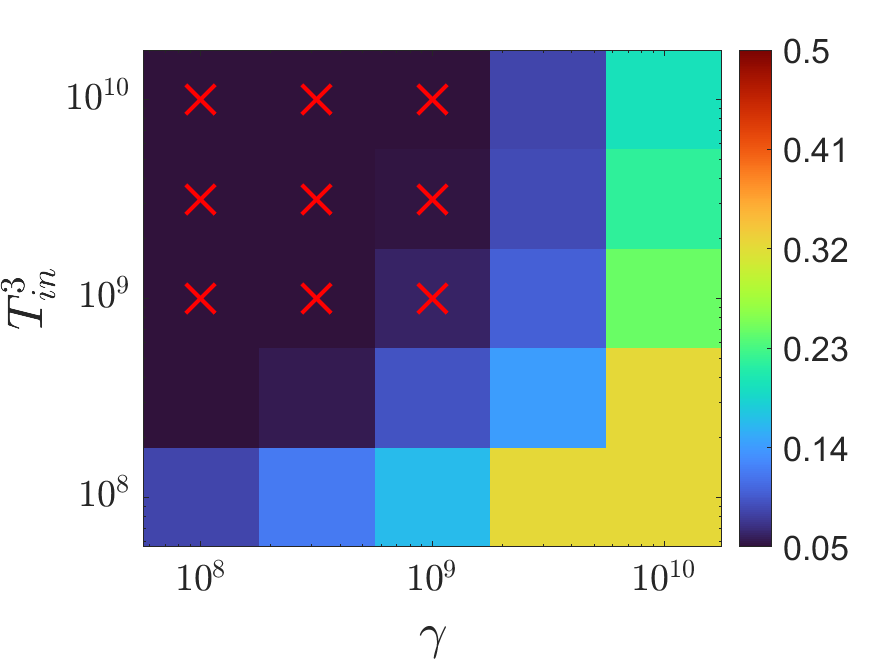}
            \caption[]%
            {{\small $\|T - \widehat{T}\|_1/\|T\|_1$}}    
            \label{fig:mean and std of net34}
        \end{subfigure}
        \begin{subfigure}[b]{0.49\textwidth}   
            \centering 
                \includegraphics[trim={0 0 0 0},clip,width=1\textwidth]{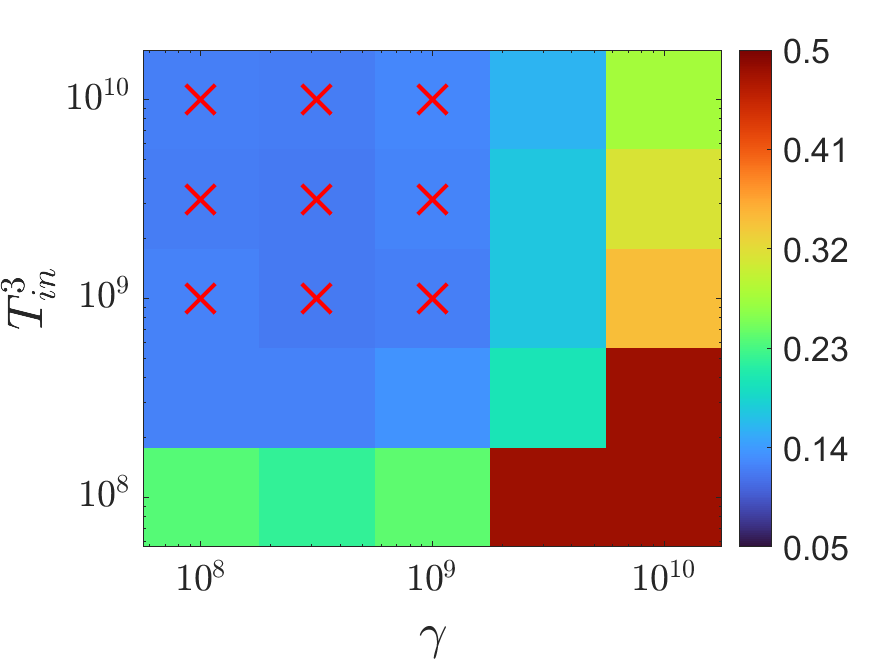}
            \caption[]%
            {{\small $\|\sigma_E E - \widehat{\sigma_E E}\|_1/\|\sigma_E E\|_1$}}    
            \label{fig:mean and std of net44}
        \end{subfigure}


\end{minipage}
&
\begin{minipage}{0.4\textwidth}

        \begin{subfigure}[b]{1\textwidth}   
            \centering 
    \includegraphics[trim={0 0 0 0},clip,width=1\textwidth]{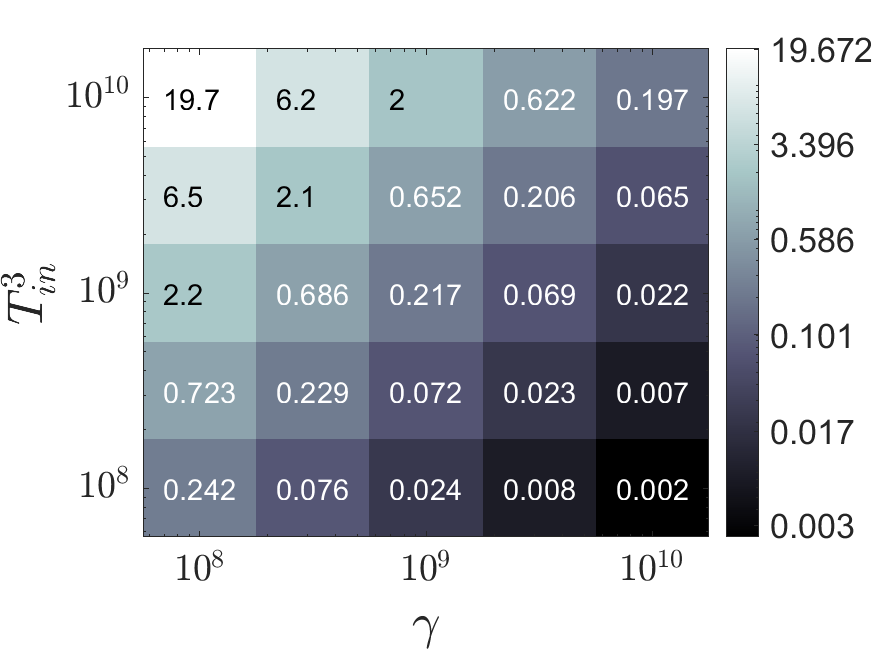}  
            \caption[]%
            {{\small $\kappa_L$}}    
            \label{fig:mean and std of net44}
        \end{subfigure}

\end{minipage}
\end{tabular}
\caption{Left: relative $L_1$ errors  ($\text{err}_{L_1}$, see eq.\ \eqref{eq:metrics}) between observed $(e,F,T,\sigma_E E)$ and learned $(\widehat{e},\widehat{F},\widehat{T},\widehat{\sigma_E E})$ solution fields over the spacetime domain $(x,t)\in [0,4\text{cm}]\times [0,2$e-$10\text{s}]$ for parameters $(\gamma, T_{in}^3) \in \{10^8,10^{8.5},10^9,10^{9.5},10^{10}\}^2$. Simulations used for training are indicated with red Xs. Right: $\kappa_L$ over the parameter range.}
\label{fig:extrap}
\end{figure}

We aim to identify to what degree the WSINDy closure model learned from the optically thin region generalizes in parameter space. It turns out that the log-linear dependence of coefficients $\Wbf$ on $(\gamma,T_{in}^3)$ identified above for model \eqref{eq:extrap_closures} enables the closure to be extrapolated in parameter space into the optically thick regime to an extent quantifiable using our dimensionless number $\kappa_L$.  

We perform this extrapolation and test on out-of-sample parameters, presenting the spacetime $L_1$ errors of the resulting closure solutions relative to the high-fidelity data ($\text{err}_{L_1}$) in Figure \ref{fig:extrap}. Plots on the left show closure errors over the $5\times 5$ parameter grid for each state variable, while the right depicts corresponding $\kappa_L$ values. To interpret the error plots, it is necessary to take into account ray effects: in the top-left region of $(\gamma,T_{in}^3)$-space, the dynamics are more optically thin and hence ray effects are more prominent. In this case, the majority of errors are due to the closure failing to capture the ray effects, which are not physical. (Here the learned closure solutions for $e$ and $F$ are actually more accurate than the high fidelity data according to $\text{err}_{L_1,j}$, see Section \ref{sec:res_effects}.) In contrast, the bottom right region is optically thick and less corrupted by ray effects, thus the error measured in this region is a more faithful measure of the fidelity of the closure. This is visualized in the solution profiles depicted in Appendix \ref{app:additionalplots}, Figures \ref{fig:extrap_worstcase1}-\ref{fig:extrap_worstcase4}. 

Taking into account these minor distortions of the $\text{err}_{L_1}$ metric due to ray effects, Figure \ref{fig:extrap} indicates that performance of the closure correlates with $\kappa_L$, and is not purely explained by distance away from the training points. The closure is most accurate over the nine training simulations, often capturing the temperature profile to within $5\%$ of the data in the $L_1$ norm, and this performance can be extended along constant $\kappa_L$ lines, and extrapolated to $\kappa_L\approx 0.07$, an unobserved parameter regime. This trend in $\kappa_L$ is most obvious for $(e,F)$, while closure solutions for $(T,\sigma_E E)$ deviate from it slightly due to underresolved inflow dynamics near the boundary producing larger errors at $(\gamma,T_{in}^3) = (10^{10},10^{10})$ (see Appendix \ref{app:additionalplots}, Figure \ref{fig:extrap_worstcase3}). Altogether, the results in Figure \ref{fig:extrap} can be used to quantify the applicable regime of the closure as $\kappa_L\gtrsim 0.07$, implying that extrapolation of the learned closures is possible within the optically thin region, but quickly breaks down in optically thick regions. 

From this we identify the six values in the bottom right corner of the parameter grid as parameters where the closure is not applicable, as here $\kappa_L < 0.07$ (see Figure \ref{fig:extrap}, right). We expect the closure to fail here, as the problem becomes optically thick because the photon mean free path is much smaller than the characteristic gradient length scale. We can draw analogies with the Knudsen number \text{kn}, where $\text{kn}\lesssim 10^{-2}$ is known to indicate fluid-like (collisional) behavior in plasma physics \cite{kramer2020plasma,vasey2025influence}. As the dynamics become more optically thick (decreasing $\kappa_L$), we start to see an overestimated flux profile (Fig.\ \ref{fig:extrap_worstcase2}), which is to be expected, as we are applying an optically thin closure to a thick region where flux is limited by diffusive behavior. When applied to the most optically thick region, the closure is unstable and blows up due to lack of positivity in $\sigma_E E$ (Fig.\ \ref{fig:extrap_worstcase4}), leading to negative temperature values (Fig.\ \ref{fig:extrap_worstcase3}). In this region a transition to known diffusive closures is expected, which we plan to explore in future work.

\subsection{Resolution effects}\label{sec:res_effects}


\begin{figure}
\begin{tabular}{c}
\includegraphics[trim={100 0 90 0},clip,width=1\textwidth]{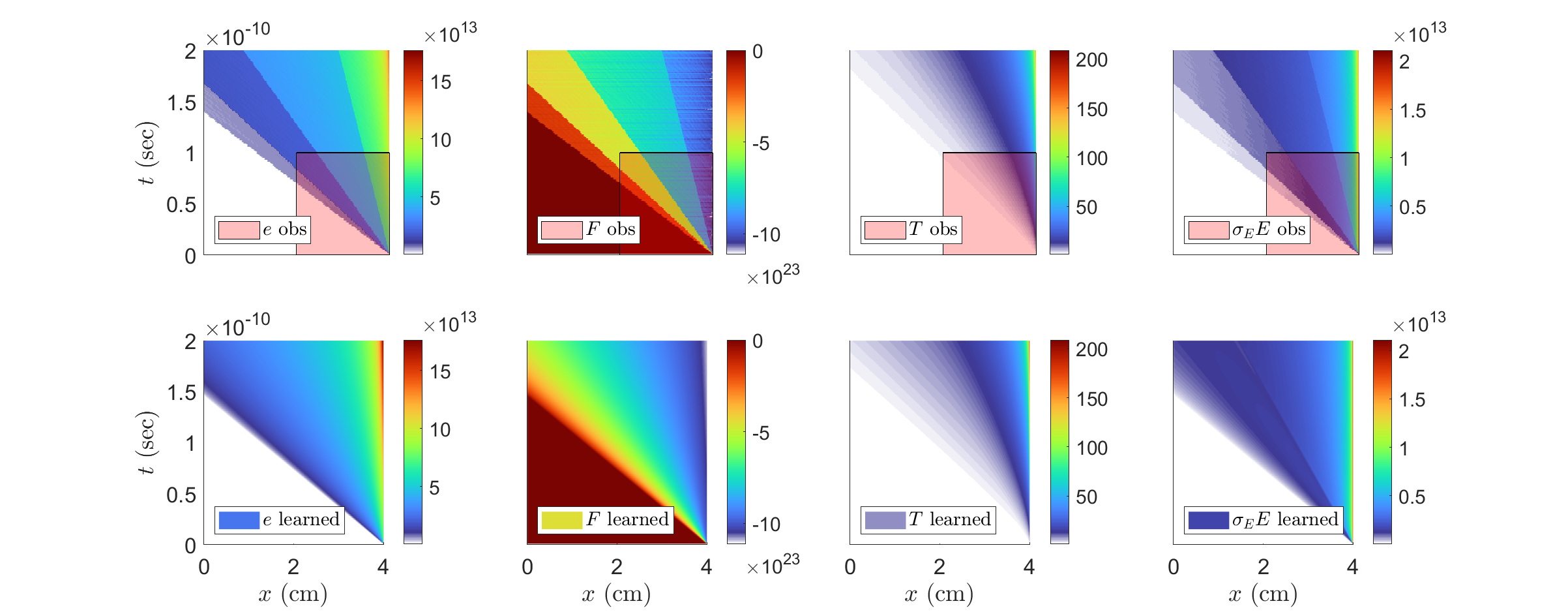} \\ \hline 
\includegraphics[trim={100 0 90 0},clip,width=1\textwidth]{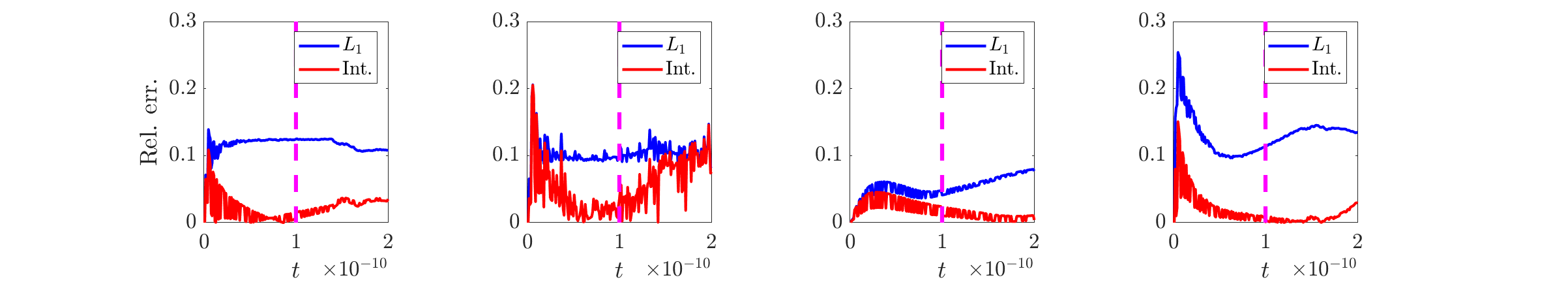} \\ 
\end{tabular}
\caption{Comparison between low-resolution high-fidelity data ($M_\Omega=8$, top row), with region in spacetime used for training boxed, and resulting closure solution (middle row) at parameters $(\gamma,T_{in}^3) = (10^9,10^9)$. The bottom row shows the relative spatial error as a function of time using the $\text{err}_{L_1,j}$ and $\text{err}_{Int,j}$ metrics (see eq. \eqref{eq:metrics}).}
\label{fig:extrap_33}
\end{figure}

\begin{figure}
\centering
\begin{tabular}{c}
\includegraphics[trim={100 0 90 0},clip,width=1\textwidth]{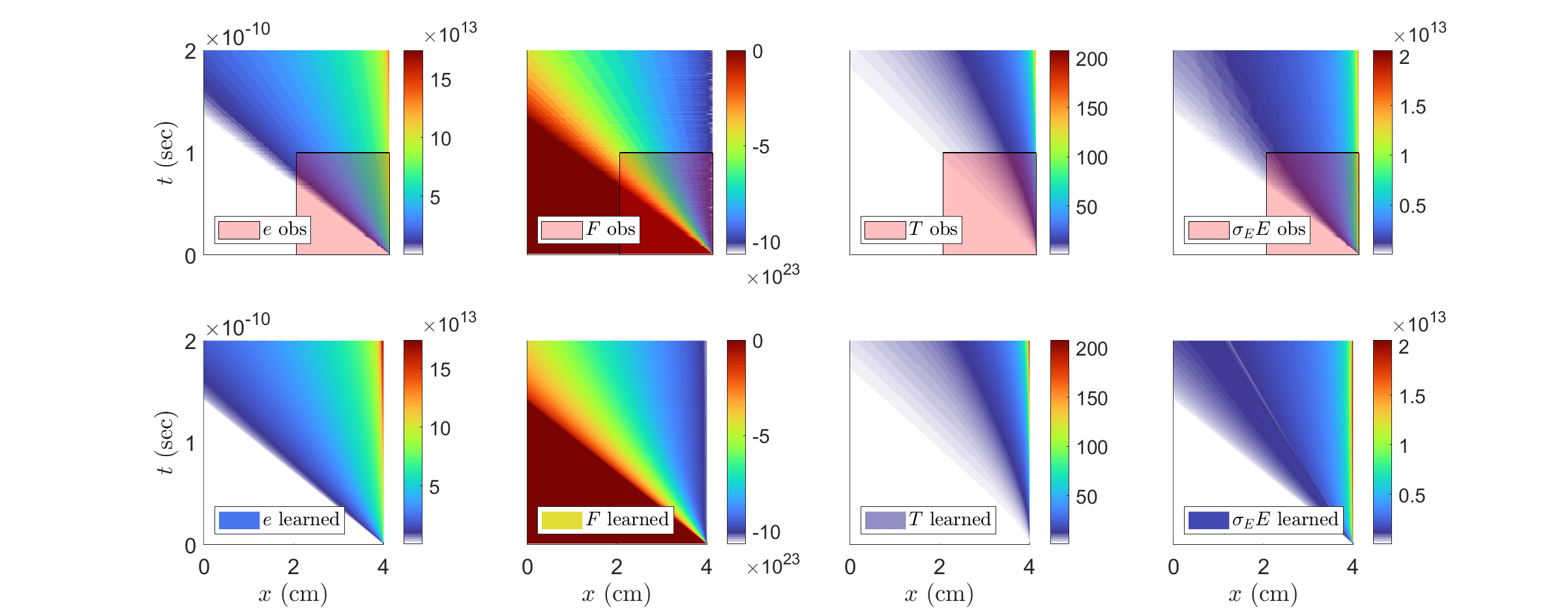} \\ \hline
\includegraphics[trim={100 0 90 0},clip,width=1\textwidth]{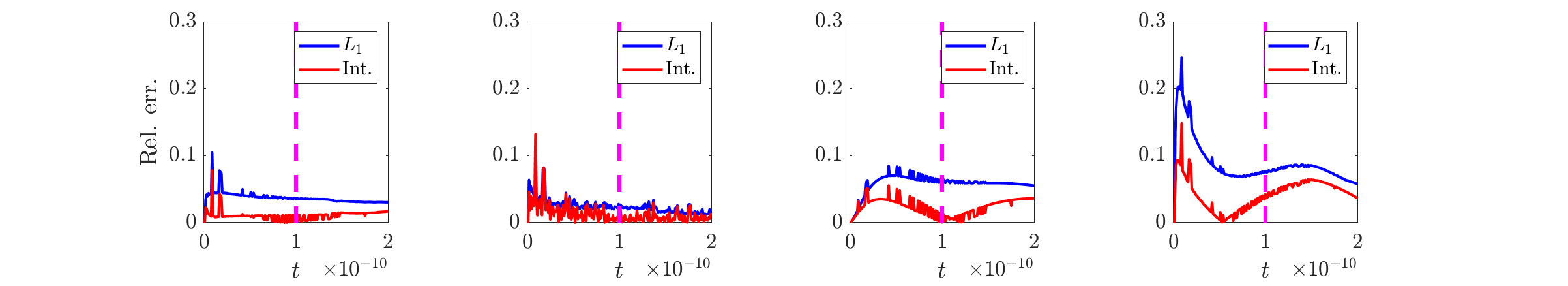} \\ 
\end{tabular}
\caption{Comparison between high-resolution high-fidelity data ($M_\Omega=48$) and resulting closure solution at $(\gamma,T_{in}^3) = (10^9,10^9)$ (compare to Figure \ref{fig:extrap_33}).}
\label{fig:highres_closure}
\end{figure}


\begin{figure}

\begin{tabular}{@{}c@{}c@{}c@{}c@{}}
\includegraphics[trim={23 0 0 0},clip,width=0.25\textwidth]{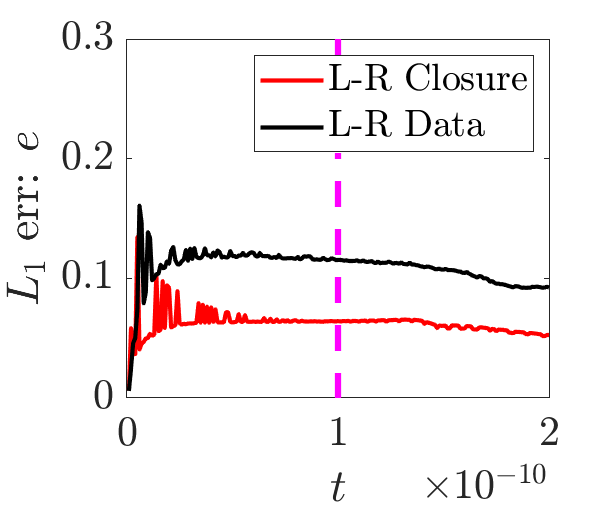}
&
\includegraphics[trim={23 0 0 0},clip,width=0.25\textwidth]{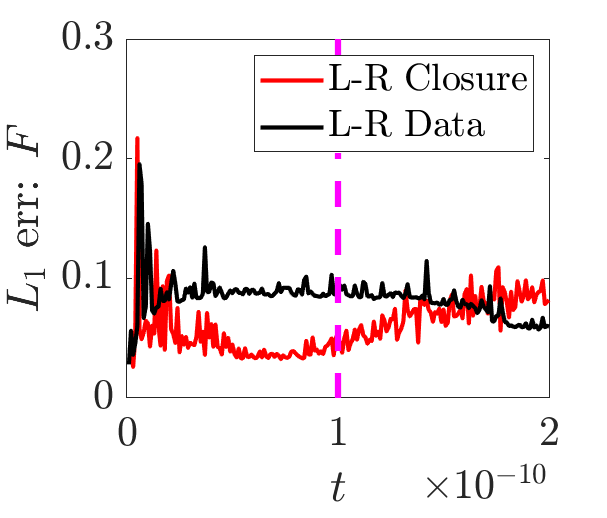}
&
\includegraphics[trim={23 0 0 0},clip,width=0.25\textwidth]{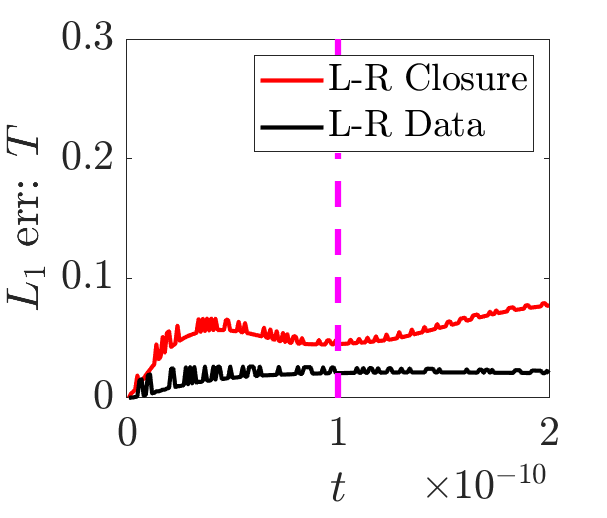}
&
\includegraphics[trim={23 0 0 0},clip,width=0.25\textwidth]{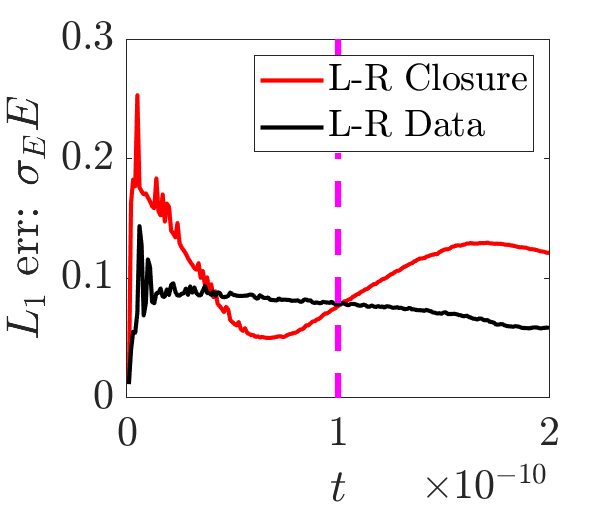} \\
\includegraphics[trim={23 0 0 0},clip,width=0.25\textwidth]{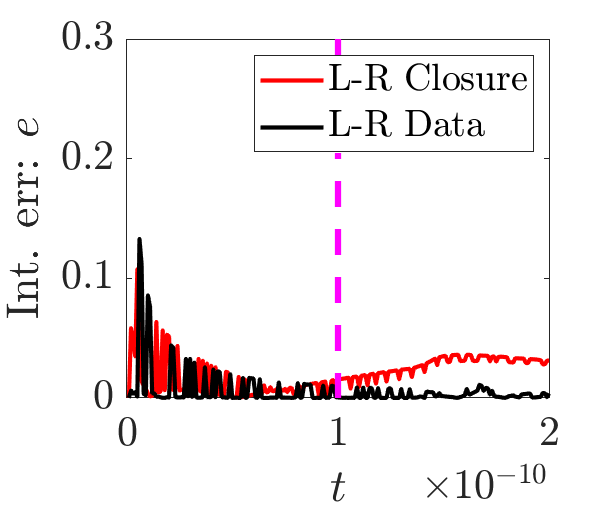}
&
\includegraphics[trim={23 0 0 0},clip,width=0.25\textwidth]{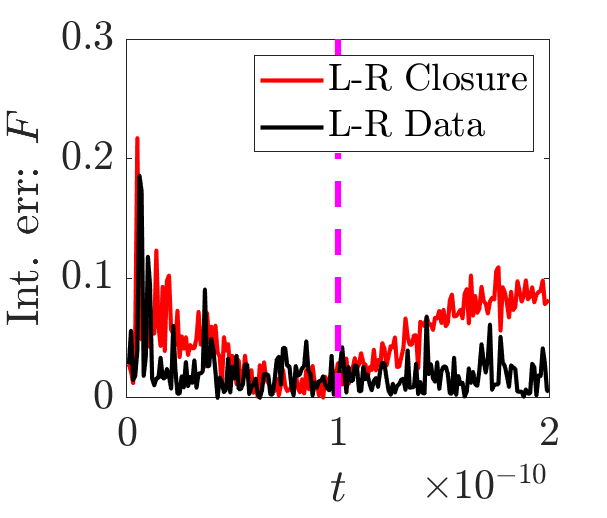}
&
\includegraphics[trim={23 0 0 0},clip,width=0.25\textwidth]{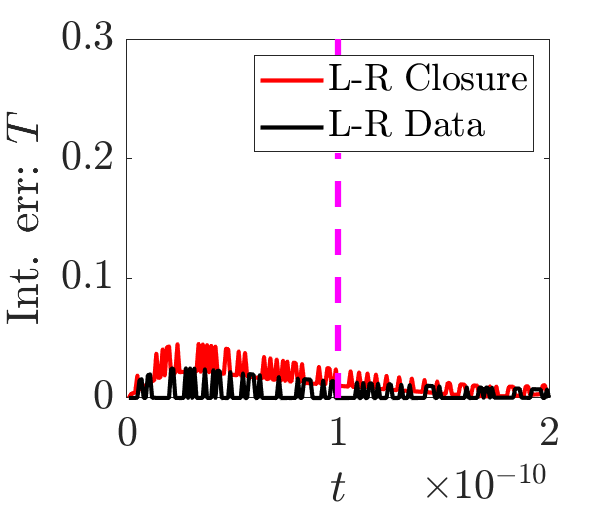}
&
\includegraphics[trim={23
0 0 0},clip,width=0.25\textwidth]{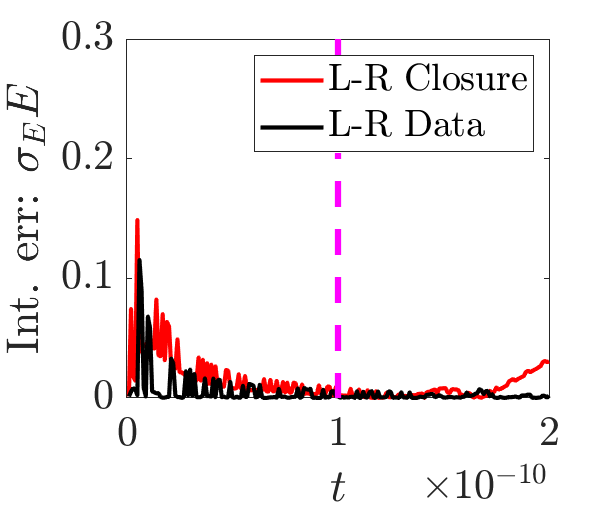} \\
$e$ & $F$ & $T$ & $\sigma_E E$
\end{tabular}

\caption{Comparison between errors incurred by the low-resolution closure solution (L-R Closure) and those of the low-resolution high-fidelity data (L-R Data, $M_\Omega=8$), both measured against the high-resolution high-fidelity data ($M_\Omega=48$) in the $\text{err}_{L_1,j}$ metric (top row) and $\text{err}_{Int,j}$ metric (bottom row). (The L-R Data was used to train the L-R closure.)}
\label{fig:high_low_res2}
\end{figure}

A surprising benefit of the proposed weak-form-based methodology is that numerical issues such as particle noise and ray effects have minor effects on the resulting closure. In fact, we find that the learned closure can be used to provide an approximation to high-resolution moment quantities that is comparable to or better than that of the lower resolution quantities, {\it even if the latter are used to train the closure model}. In this way, optically thin WSINDy closures for TRT have the potential to increase the accuracy of high fidelity simulations. 

To study this, we compare the use of low-resolution ($M_\Omega=8$) training data as above to high-resolution training data, defined as high fidelity simulations with $M_\Omega=48$ polar angles. As a reference point, Figure \ref{fig:extrap_33} depicts the performance of the closure model (middle row) in comparison with the low-res data (top row) at parameters $(\gamma,T_{in}^3) = (10^9,10^9)$. The bottom row shows the relative $L_1$ error ($\text{err}_{L_1,j}$) and the relative integral error ($\text{err}_{Int,j}$) as functions of time.  Ray effects manifest as a significant difference between the two error metrics, with the relative $L_1$ error on the order of $10\%$ over the course of the time series (slightly lower for $T$), while the integral error is closer to $2\%$. Particle noise contributes an additional portion to both error quantities, most notably in $F$. (Note that the largest errors occur at early times when the solution has not propagated far into the domain, i.e. the early-time errors arise from a small denominator). 

In Figure \ref{fig:highres_closure}, we instead learn the closure from high-resolution high-fidelity data at the same parameter values. Here we see that the relative $L_1$ error decreases significantly, as ray effects present smaller (yet more frequent) solution jumps. The overall error between the closure solution and training data clearly decreases with resolution, but more importantly, {\it the closure solutions remain nearly identical} (as can be seen by comparing the middle rows of \ref{fig:extrap_33} and \ref{fig:highres_closure}), while the low-res and high-res high-fidelity solutions are quite different (recall Fig.\ \ref{fig:ray_effects}). This indicates that the low-res closure solution is closer to the high-res high-fidelity solution than the low-res one. While this is qualitatively clear from solution plots in Figures \ref{fig:extrap_33} and \ref{fig:highres_closure}, we quantify this in Figure \ref{fig:high_low_res2}, which compares the relative error between the high-res high-fidelity solution and the closure solution (trained with the low-res high-fidelity data, red curves), with the relative error between the low and high resolutions datasets (black curves), in the $L_1$ and integral metrics (top and bottom rows). It is clear from the bottom row that in the $\text{err}_{Int,j}$ metric, the low-res data approximates the high-res data slightly better than the closure solution does, although all such errors are reasonably low. In the $L_1$ metric (top row), however, the closure produces a more accurate solution to $e$ and $F$, with errors in $e$ approximately half that of the low-res data vs.\ the high-res data. Further, the $L_1$ errors of the closure solution against the high-res data are lower than the errors against the low-res data (Fig.\ \ref{fig:extrap_33}). This indicates that the solutions from the closure model trained on low-res high-fidelity data are actually closer to the high-res high-fidelity data than the low-res data. The closure solutions for $T$ and $\sigma_E E$ (rightmost two columns of Fig.\ \ref{fig:high_low_res2}) exhibit larger errors than that of the low-res data, indicating that improvements to the closure for $\sigma_E E$ are possible. Judging by the error decrease in $\sigma_E E$ going from low-res to high-res data (Figs.\ \ref{fig:extrap_33}, \ref{fig:highres_closure}), we expect this to be overcome by increasing the fidelity and amount of data, as well as the size of libraries employed for regression. 

\section{Discussion}

We have demonstrated that WSINDy is a viable framework for finding interpretable data-driven closures for thermal radiation transport. We have chosen to focus on the notoriously challenging regime of optically thin media, in which the widely-used suite of diffusive closures is invalid. We have further restricted examples to a set of multifrequency radiation problems to demonstrate that, by treating the energy-weighted opacity $\sigma_E E$ as a state variable, we arrive at high-fidelity closures that correctly account for this evolving quantity which differs greatly from its gray limit $\propto T^{-3}$. WSINDy readily finds an evolution equation for $\sigma_E E$ indicating that this nonlocal quantity can be treated locally by augmenting the state space.

Using a simple change of variables to total energy $e$, we are able to derive libraries and constraints that promote model stability and physical properties, namely hyperbolicity, source stability, black-body equilibria, and rotational symmetry. These appear as linear constraints on the WSINDy coefficients or library constraints, the latter being automatically enforced in any model, and do not considerably increase the runtime of WSINDy (which runs in seconds to minutes on a laptop). We anticipate that further model constraints may be enforced in a similar manner. 

A crucial element of our closure learning method is the weak form, in particular due to the quality of high-fidelity data employed for training. We predominantly discretize angle space with $M_\Omega=8$ discrete polar angles, which leads to severe ray effects that worsen in the optically thin regime, evidence that the underlying distribution is highly anisotropic in angle space. In addition, we are affected by deterministic particle discretization errors, a similar but less severe corruption than statistical noise. It has been shown previously that WSINDy is capable of identifying mean-field equations from particle data \cite{messenger2022learning}, and is effective at identifying PDEs from non-smooth data \cite{messenger2020weakpde}.  However, the current article presents a different use-case: identifying closures that {\it remove} these numerical artifacts without attempting to model them. The effective fidelity of the closures presented appears to be closer to the $M_\Omega=48$ polar angle discretization than the $M_\Omega=8$ used in training (Figures \ref{fig:ray_effects}, \ref{fig:extrap_33}, \ref{fig:high_low_res2}). 
A goal for future work is to combine our WSINDy-based data-driven closures with HOLO iterations to further accelerate solutions to TRT, in combination with reducing numerical artifacts such as ray effects. Overall, we aim to investigate the utility of the presented results in accelerating high-fidelity simulations. 

Compared to black-box machine learning methods, we have demonstrated the flexibility and versatility of weak-form, library-based closure learning in enforcing constraints and dealing with challenging non-Gaussian corruptions, opening the door to other kinetic closure problems. We have also demonstrated that such closures can be {\it extrapolated} in ways that are not common with black-box methods. We have extrapolated our closures in space and time, training on data from a restricted space-time domain and simulating the learned model over a larger spatial domain for a longer time window. This is not expected to be easily achieved by black-box representations. We have also demonstrated that our closures can be extrapolated in {\it parameter space}, and provided explicit coefficient parameter relations in Table \ref{tab:coeffs}. These are not surprising, and in hindsight are intuitive given the dependence on $\kappa_L$, the dimensionless number we have implicated in defining the relevant region of parameter space for the closures presented. Through these relations, we are able to establish regions in parameter space where the proposed closure is valid, and regions where it is not valid. While the utility of black-box methods is unquestionable, we emphasize the power of interpretable methods such as those discussed here in aiding in understanding, remaining computationally efficient, and offering a path to generalization.  

We acknowledge that the examples considered are far from a complete treatment of the TRT dynamics relevant to engineering applications, especially when coupled with additional physics such as hydrodynamics. We have kept the exposition simple in order to demonstrate a targeted use-case for WSINDy-based closures that preserve physical properties: optically thin, multifrequency dynamics in 1D, a regime without a known effective moment closure approximation. We anticipate that extensions to multiple dimensions and multiple materials will bring challenges of their own, which we look forward to tackling in future work. We have also not solved the boundary closure problem, which arises because boundary data at the moment level in general depends on the kinetic solution. For the problems considered here, boundary data can be extrapolated from the training simulation boundary data in a similar manner as the closure parameters (not shown\footnote{The training boundary data is already fit to a quadratic function of time to eliminate particle noise in the examples above. The parameters of the quadratic can be extrapolated over $(\gamma,T_{in}^3)$.}). For other TRT problems, it will be necessary to find a boundary closure. We leave this to future work as well. Finally, we acknowledge that in the limited setting considered here, our closures do still accrue errors when compared to high fidelity data, especially in $\sigma_E E$. Since moment closure is an inexact science, we expect errors to result from the averaging of microscale quantities. We also note that our technique of learning an auxiliary evolution equation for an opacity moment appears to be novel (to the best of the authors' knowledge, it is not found in the literature), so we expect  improvements on it to be made in future work. We have presented evidence that our results improve with better quality data, and we expect the same with library refinement and other improvements in optimization. From the starting point presented here, we anticipate much in the form of improved kinetic closures and successful combinations with forward solvers.   

\section*{Acknowledgments}
  This work was supported by the Laboratory Directed Research and Development program and the U.S. Department of Energy Office of Advanced Scientific Computing
  Research through the CHaRMNET Mathematics Multifaceted Integrated Capability Center and the ASCR Competitive Portfolios program at Los Alamos National
  Laboratory (LANL), under contract No. 89233218CNA000001.  
\appendix

\section{Closure solution plots}\label{app:additionalplots}
\begin{figure}
\begin{center}
\includegraphics[trim={100 10 100 30},clip,width=1\textwidth]{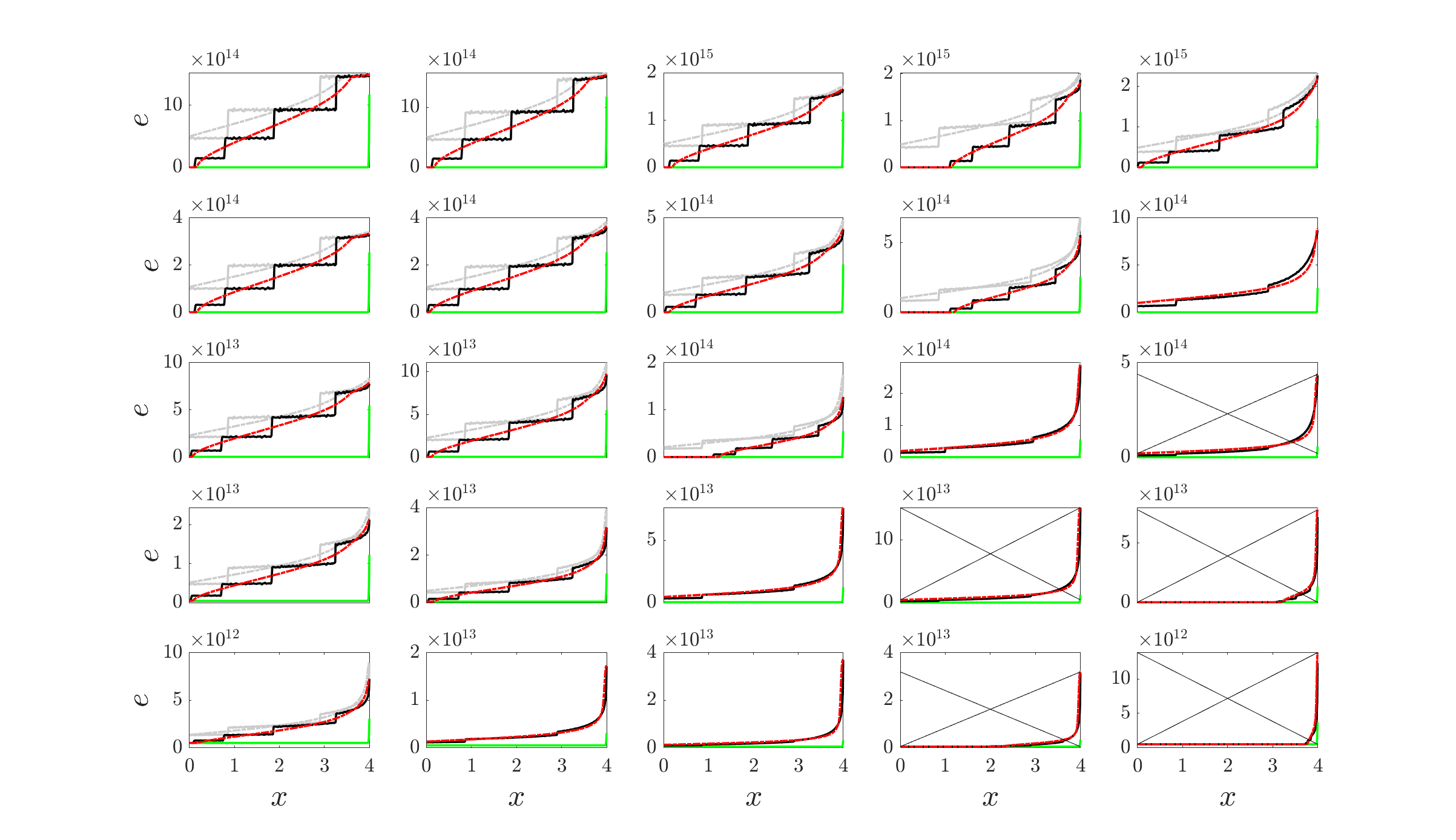}
\end{center}
\caption{Plots of closure solution for total energy $e$ (red) vs high-fidelity data (black) that achieves the worst-case relative spatial $L_1$ error ($\text{err}_{L_1,j}$, see eq. \eqref{eq:metrics}) over the extrapolation region $t\in[1$e-$10,2$e-$10]$ of the time domain. Plots correspond with $(\gamma,T_{in}^3)$ values consistent with Figure \ref{fig:extrap}. For reference, initial conditions are shown in green, while closure solutions and data at the final time ($t=2$e-$10$) are shown in gray (overlapped by the worst-case plots in several cases). Plots with X's are failure modes where the closure generated negative $\sigma_E E$ values.}
\label{fig:extrap_worstcase1}
\end{figure}

\begin{figure}
\begin{center}
\includegraphics[trim={95 10 100 30},clip,width=1\textwidth]{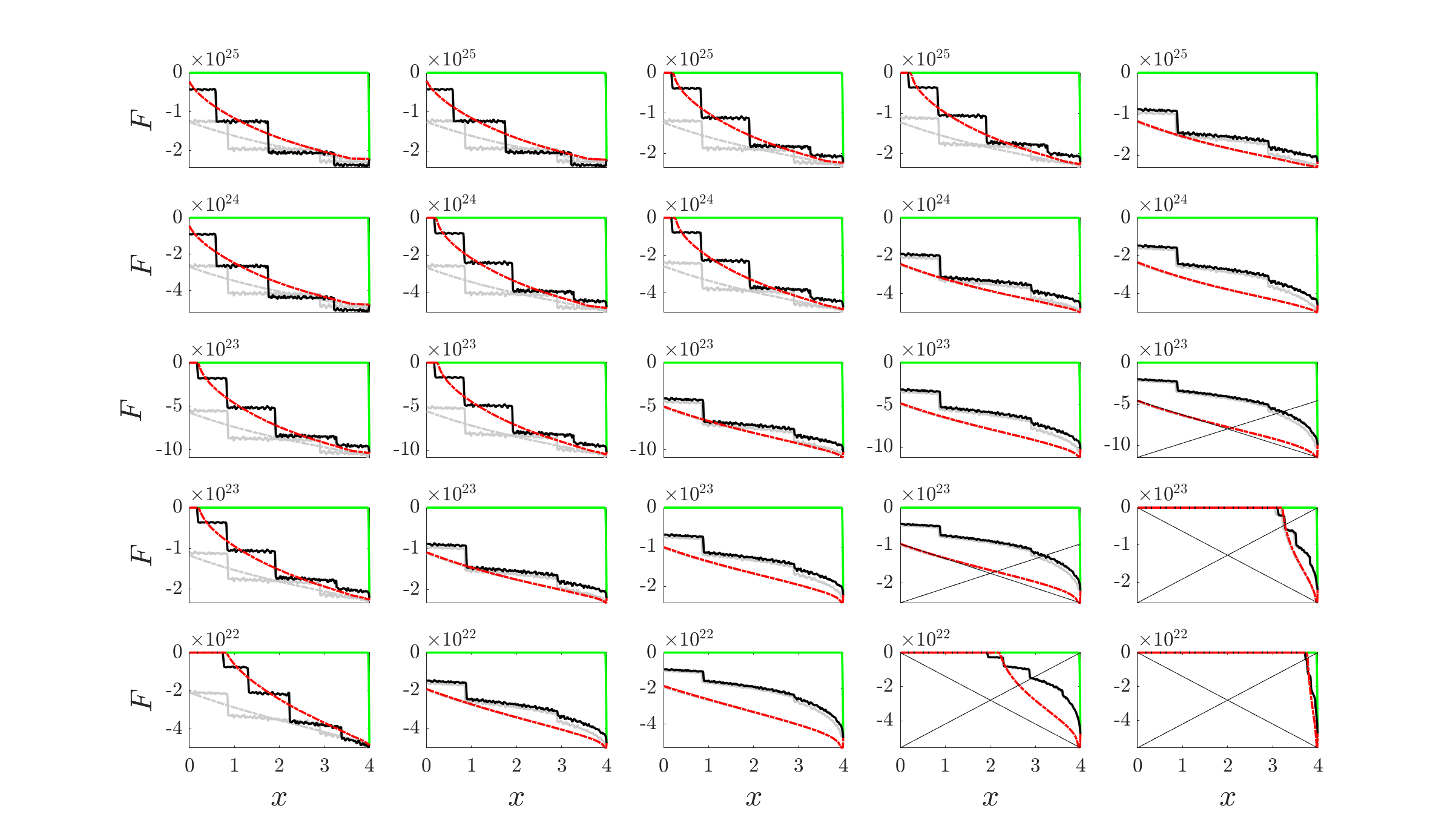}
\end{center}
\caption{Plots of closure solution for radiation flux $F$ (red) vs data (black) that achieves the worst-case relative spatial $L_1$ error over the extrapolation region in time (analogous to Figure \ref{fig:extrap_worstcase1}).}
\label{fig:extrap_worstcase2}
\end{figure}

\begin{figure}
\begin{center}
\includegraphics[trim={90 10 100 30},clip,width=1\textwidth]{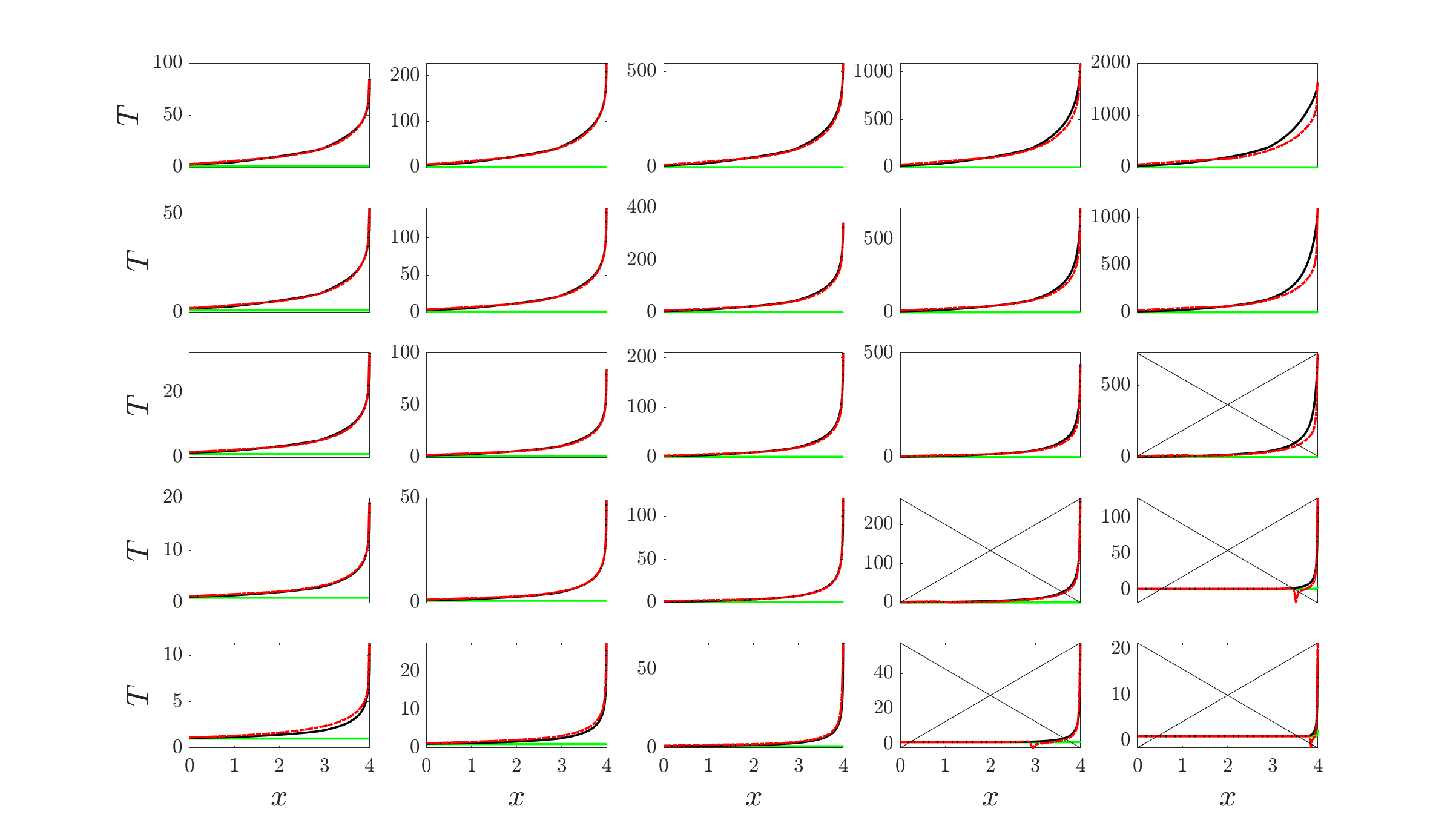}
\end{center}
\caption{Plots of closure solution for temperature $T$ (red) vs data (black) that achieves the worst-case relative spatial $L_1$ error over the extrapolation region in time (analogous to Figure \ref{fig:extrap_worstcase1}).}
\label{fig:extrap_worstcase3}
\end{figure}

\begin{figure}
\begin{center}
\includegraphics[trim={100 10 100 30},clip,width=1\textwidth]{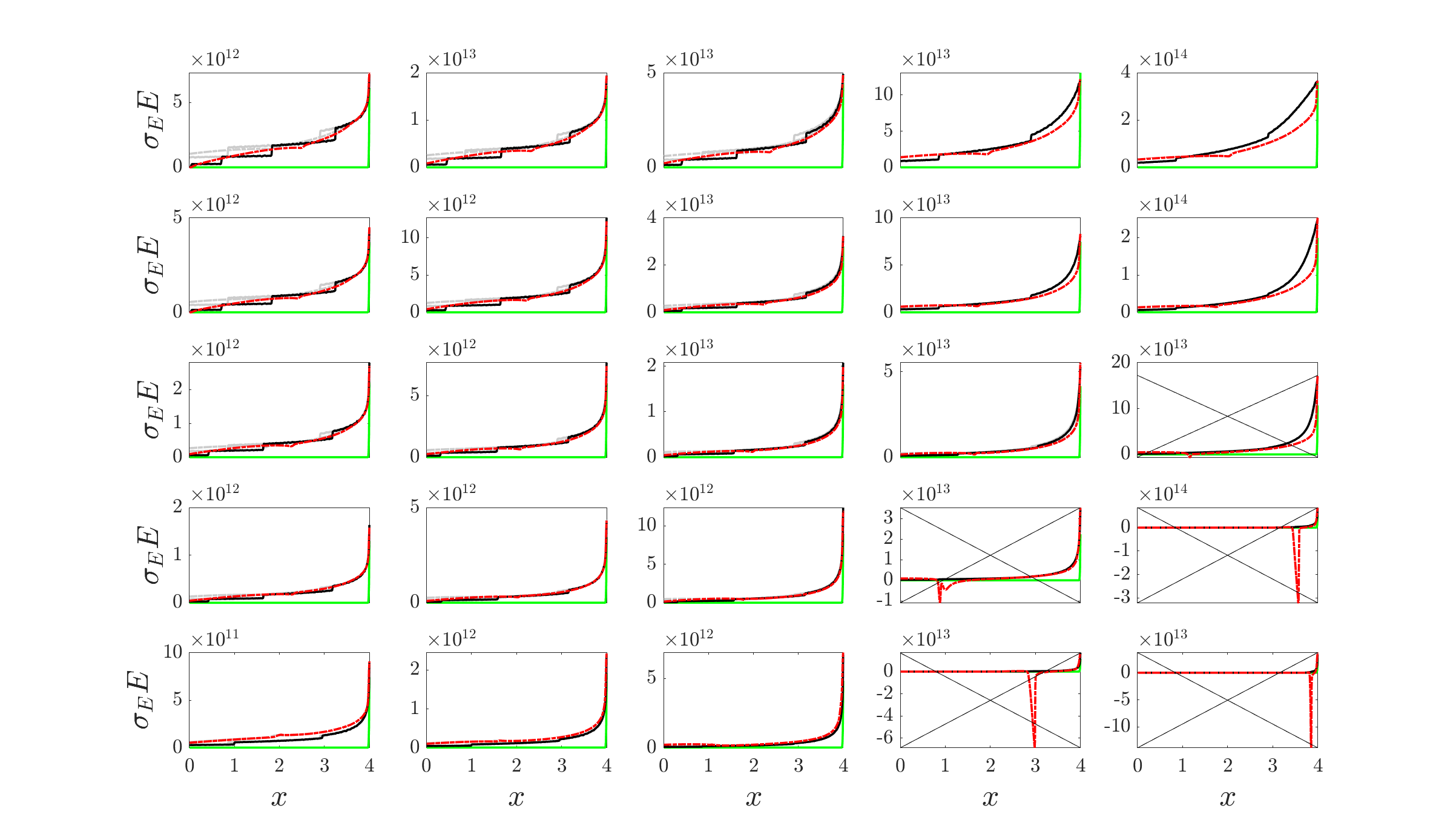}
\end{center}
\caption{Plots of closure solution for energy-weight opacity $\sigma_E E$ (red) vs data (black) that achieves the worst-case relative spatial $L_1$ error over the extrapolation region in time (analogous to Figure \ref{fig:extrap_worstcase1}).}
\label{fig:extrap_worstcase4}
\end{figure}

Figures \ref{fig:extrap_worstcase1}-\ref{fig:extrap_worstcase4} contain depictions of the worst-case-in-time relative spatial $L_1$ errors of solutions to the closure model \eqref{eq:extrap_closures}. From Figures \ref{fig:extrap_worstcase1}, \ref{fig:extrap_worstcase2}, it can be seen that ray effects of the $M_\Omega=8$ high-fidelity data (black) represent significant corruptions, which do not appear in the closure solution (red). It is also clear in each figure that that similar solution behavior is observed along $\kappa_L=\text{const.}$ lines (diagonals). 

\section{Review of known TRT closures}\label{sec:known_closures}

\bi{Isotropic intensity:} If radiation intensity $I$ does not depend on angle $\bOmega$, then the Eddington tensor has a simple closure, $\CalE = \frac{c}{3}\Ibf_d$ where $d$ is the spatial dimension and $\Ibf_d$ is the identity in $\Rbb^d$. This leads to the $P_1$ system:
\begin{align}
    \frac{\partial E}{\partial t} &= -\nabla \cdot \mathbf{F} - c\sigma_E(T)E + ac\sigma_P(T) T^4,\\
    \frac{1}{c}\frac{\partial \mathbf{F}}{\partial t} & = -\frac{c}{3}\nabla E - \sigma_R(T) \mathbf F , \label{eq:moment-eqs-F-uni}\\
    \rho c_V \frac{\partial T}{\partial t} & = -ac\sigma_p(T)T^4 + c\sigma_E(T) E.
\end{align}
This is also obtained from a first-order Legendre expansion of $I$. In the optically thin (streaming) case, the wavespeed is underestimated at $c/\sqrt{3}$. This is corrected in the $P_{1/3}$ system \cite{olson2000diffusion}. Typically these are applied with the additional assumption of gray transport.\\

\noindent \bi{Gray transport:} Gray radiation refers to eliminating frequency from phase space. There are several ways to do this. For simplicity we present here a ``single-frequency'' model. If all photons have the same frequency $\nu_0$, then $I = I_0(\xbf,\bOmega,t) \delta(\nu-\nu_0)$. Integrating the TRT equations and defining $\sigma_0(T) := \sigma(\nu_0,T)$, we get
\begin{subequations}\label{eq:ho_gray}
\begin{align}\label{eq:ho-I}
  \frac{1}{c} \frac{\partial I_0}{\partial t} & = - \bOmega \cdot \nabla I_0 - \sigma_0(T) I_0 + \tfrac{1}{4 \pi}ac\sigma_0(T)T^4,\\
  \rho c_V \frac{\partial T}{\partial t} & = \sigma_0(T)\left(\int_{S^d} I_0\d\bOmega\right)  - ac\sigma_0(T) T^4.
\end{align}
\end{subequations}

\noindent \bi{Gray Radiation Diffusion:} The Gray TRT equations are still nonlocal, hence benefit from a moment expansion into $E,\Fbf,T$. If we assume that $\partial_t \Fbf \ll c$, then solve for $\Fbf$ in \eqref{eq:lo} gives \cite{levermore1981flux}
\[\Fbf \approx -\frac{1}{\sigma_0(T)}\CalE \nabla E\]
and the moment system reduces to
\begin{subequations}
\begin{align}
    \frac{\partial E}{\partial t} &= \nabla \cdot \left(\frac{1}{\sigma_0(T)}\CalE \nabla E\right) - c\sigma_0(T) E + ac\sigma_0(T) T^4,\\
    \rho c_V \frac{\partial T}{\partial t} & = c\sigma_0(T) E-ac\sigma_0(T)T^4.
\end{align}
\end{subequations}
\noindent \bi{Gray-Radiation, Symmetric-Diffusion:} 
Things simplify if we further assume that variations in the direction of travel are isotropic, or for some function $v$ corresponding to the scaled angular variance, or
\[v^2(\xbf,t) = \CalE_{ii} = c\frac{\int_{S^d}(\bOmega_i)^2I_0 d\bOmega}{\int_{S^d}I_0d\bOmega} = \frac{\int_{S^d}(\bOmega_i)^2I_0 d\bOmega}{E},\]
we have $\CalE_{ij}=0$ and so
\[\CalE = v^2(\xbf,t)\Ibf_d \]
Then 
\[\frac{\partial E}{\partial t} = \nabla \cdot \left(\frac{v^2}{\sigma_0(T)}\nabla E\right) - c\sigma_0(T) E + ac\sigma_0(T) T^4\]
To the best of the author's knowledge this case has not been considered in the literature, but we include it here as another example of an analytic closure. A closed equation for $v^2$ can be found if we assume that third angular moments of $I_0$ are spatially constant. We have
\begin{align}
  \partial_t\int\Omega_i^2I_0d\bOmega = c\partial_t(Ev^2) &= c(\partial_t E)v^2 + cE\partial_t(v^2) \\
                    &\approx c\left(\nabla\cdot (\frac{v^2}{\sigma_0}\nabla E) - c\sigma_0 E+ ac \sigma_0 T^4\right)v^2 +cE\partial_t (v^2) 
                    \intertext{but also, using the 3rd moment property,}
\partial_t\int\Omega_i^2I_0d\bOmega &= c\int\Omega_i^2\left(-\bOmega\cdot\nabla I_0 - \sigma_0I_0 +\frac{ac}{\omega_d} \sigma_0T^4\right)d\bOmega \\
                    &= -c^2\sigma_0 Ev^2 + \frac{ac^2}{3}\sigma_0 T^4 
\end{align}
Equating the two and dividing by $c Ev^2$ gives
\[\partial_t y = -\frac{1}{E}\nabla \cdot \left(\frac{e^y}{\sigma_0}\nabla E\right) -\frac{ac\sigma_0 T^4}{E}\left(1-\frac{1}{3e^y}\right)\]
for $y := \log(v^2)$.\\

\noindent \bi{Gray-Radiation Isotropic-Diffusion:} 
If we assume that $v^2 = \text{const}$ above, then the system closes with the usual gray radiation diffusion equations: 
\begin{subequations}\label{eq:rd}
\begin{align}
    \frac{\partial E}{\partial t} &= \nabla \cdot \left(\frac{v}{\sigma_0(T)}\nabla E\right) - c\sigma_0(T) E + ac\sigma_0(T) T^4,\\
    \rho c_V \frac{\partial T}{\partial t} & = c\sigma_0(T) E-ac\sigma_0(T)T^4.
\end{align}
\end{subequations}
and typically $v = c/3$ (the approximation made when $I$ is angularly independent).\\

\section{Overview of WSINDy}\label{app:WSINDy_overview}
Assume that a desired PDE for $u$ takes the form
\begin{equation}\label{pdemodel}
D^{\pmb{\alpha}^{(0)}}u(\xbf,t) = \sum_{i,j=1}^{I,J} \wstar_{(i-1)J+j} D^{\pmb{\alpha}^{(i)}}f_j(u(\xbf,t)), \quad (\xbf,t)\in \Omega\times[0,\infty),
\end{equation}
where $\Omega\subset\Rbb^d$ is a bounded open set. The operators $D^{\pmb{\alpha}^{(i)}}$ for $1\leq i \leq I$ represent any linear differential operator in the variables $(\xbf,t)\in \Rbb^{d+1}$, where $\pmb{\alpha}^{(i)}=(\pmb{\alpha}^{(i)}_1,\dots,\pmb{\alpha}^{(i)}_{d+1})$ is a multi-index such that 
\[D^{\pmb{\alpha}^{(i)}}v = \frac{\partial^{\pmb{\alpha}^{(i)}_1+\cdots+\pmb{\alpha}^{(i)}_d+\pmb{\alpha}^{(i)}_{d+1}}}{\partial \xbf_1^{\pmb{\alpha}^{(i)}_1}\cdots\partial \xbf_d^{\pmb{\alpha}^{(i)}_d}\partial t^{\pmb{\alpha}^{(i)}_{d+1}}}v.\]
The left-hand side operator $D^{\pmb{\alpha}^{(0)}}$ is assumed to be known. The functions $f_j:\Rbb\times\Rbb^d \to \Rbb$, $1\leq j\leq J$, contain a basis for nonlinearities present in the model, and together with the linear operators $D^{\pmb{\alpha}^{(i)}}$ comprise the feature library $\Theta:=\{D^{\pmb{\alpha}^{(i)}}f_j\}_{i,j=1}^{I,J}$. The weight vector $\wstar \in \Rbb^{IJ}$ is assumed to be {\it sparse} in $\Theta$. We assume that we have access to approximate solution data
\begin{equation}\label{datassump}
\Ubf(\Xbf,\tbf) = u(\Xbf,\tbf)+\ep
\end{equation}
where $u$ solves \eqref{pdemodel} for some weight vector $\wstar$, and that $(\Xbf,\tbf) \in \Rbb^{n_1\times \cdots \times n_d\times n_{d+1}}$ is a fixed known spatial-temporal grid of points in $\Omega \times [0,T)$ having $n_i$ points in the $i$th dimension. Here $\ep$ represents mean-zero noise, typically assumed i.i.d.\ with fixed variance $\sigma^2<0$ associated with sampling the underlying solution $u(x,t)$ at any point $(x,t)\in \Omega\times[0,T)$. 

The problem is then to identify the sparse weight vector $\wstar$. The WSINDy algorithm \cite{messenger2020weak,messenger2020weakpde}  proposes to solve this problem by first convolving equation \eqref{pdemodel} with a smooth test function $\psi(\xbf,t)$, compactly supported in $\Omega\times [0,T]$. After integrating by parts to put all partial derivatives onto $\psi$, this leads to the {\it convolutional weak form}:
\begin{equation}\label{weakform}
D^{\pmb{\alpha}^{(0)}}\psi*u(\xbf,t) = \sum_{i,j=1}^{I,J}\wstar_{(i-1)J+j} D^{\pmb{\alpha}^{(i)}}\psi*f_j(u)(\xbf,t),
\end{equation}
where convolutions are performed over space and time. For efficiency, the test function $\psi$ is chosen to be separable,
\begin{equation}\label{septest}
\psi(\xbf,t)=\phi_1(\xbf_1)\cdots\phi_d(\xbf_d)\phi_{d+1}(t) .
\end{equation}
with each $\phi_i$ chosen according to the Fourier spectrum of $\Ubf$ to reduce the effects of high-frequency noise as in \cite{messenger2020weakpde}. Briefly, this entails first computing a wavenumber $k_i$ that approximately separates signal- and noise-dominated modes of $\Ubf$ in each of the $i\in \{1,\dots,d+1\}$ dimensions. The user then chooses two hyperparameters $\tau$ and $\hat{\tau}$, where $\tau$ is the value of $\phi_i$ at the boundary of its support (it is assumed that $\phi_i(0)=1$ and $\phi$ is symmetric) and $\hat{\tau}$ is the number of standard deviations into the power spectrum of $\phi_i$ (normalized as a probability over wavenumbers) where $k_i$ occurs. As such, $(\tau,\hat{\tau})$ sets the decay rates of $\phi_i$ in real and Fourier space (respectively), where the former increases integration accuracy and the latter reduces the effects of high-frequency noise. In \cite{messenger2020weakpde} default values of $(\tau,\hat{\tau}) =(10^{-10},2)$ are proposed, although deviations from i.i.d. Gaussian noise and/or low-regularity (nonsmooth) data may require further tuning of these hyperparameters.

Once $\psi$ is chosen, the convolutional weak form is sampled over a finite set of {\it query points} $\CalQ:=\{(\xbf^{(q)},t_q)\}_{q=1}^Q\subset \Omega\times (0,T)$, with $u$ replaced by $\Ubf$. Convolutions can be efficiently computed using the fast Fourier transform (FFT), which, due to the compact support of $\psi$, is equivalent to the trapezoidal rule, a highly-accurate scheme on noise-free data ($\ep = 0$). This results in a linear system 
\[\bbf \approx \Gbf\wstar,\]
 where the $q$th entry of $\bbf$ is $\bbf_q = D^{\pmb{\alpha}^{(0)}}\psi*\Ubf(\xbf^{(q)},t_q)$ and $q$th entry of the $((i-1)J+j)$th column of $\Gbf$ is $\Gbf_{q,(i-1)J+j}=D^{\pmb{\alpha}^{(i)}}\psi*f_j(\Ubf)(\xbf^{(q)},t_q)$. As the library $\Theta$ grows, sparsity of the weight vector $\what$ becomes necessary in order to interpret and efficiently simulate the resulting PDE. Using the assumption that $\wstar$ is sparse, the optimization problem for $\what \approx \wstar$ becomes the sparse recovery problem
 \begin{equation}\label{sparserec}
\min_{\wbf\in \Rbb^{IJ}} F(\wbf;\lambda) = \min_{\wbf\in \Rbb^{IJ}} \frac{1}{2}\nrm{\Gbf\wbf-\bbf}_2^2 + \frac{1}{2}\lambda^2\nrm{\wbf}_0.
\end{equation}
The sparsity threshold $\lambda>0$ must be set by the user and is designed to strike a balance between fitting the data, associated with low residual $\nrm{\Gbf \wbf-\bbf}_2$, and finding a parsimonious model, indicated by low $\nrm{\wbf}_0$ (and its value is typically calibrated via cross-validation)  \cite{hastie2009elements,10.5555/2526243}. A major challenge is that the columns of $\Gbf$ are typically highly correlated since they are each constructed from the same dataset $\Ubf$, which leads to many popular algorithms for solving \eqref{sparserec} performing poorly, such as convex relaxation using the $\ell_1$-norm \cite{meinshausen2006high,fan2014endogeneity}. The following approach is taken for WSINDy, which has proved to be successful under various noise levels and systems of interest, and has a theoretical guarantee proven in \cite{messenger2024asymptotic}. For $\lambda>0$ define the inner sequential thresholding step
\begin{equation}\label{MSTLS1}
\text{MSTLS}(\Gbf,\bbf; \lambda\,)\qquad \begin{dcases} \hspace{0.3cm}\wbf^{(0)} = \Gbf^\dagger \bbf \\ \hspace{0.43cm}\CalI^{(\ell)} = \{1\leq k\leq IJ\ :\ L_k(\lambda)\leq|\wbf^{(\ell)}_k|\leq U_k(\lambda)\} \\
\wbf^{(\ell+1)} = \argmin_{\supp{\wbf}\subset \CalI^{(\ell)}} \nrm{ \Gbf  \wbf-\bbf}_2^2.\end{dcases}
\end{equation}
Letting $\Gbf_{k}$ be the $k$th column of $\Gbf$, the lower and upper bounds are defined
\begin{equation}\label{MSTLSbnds} 
\begin{dcases} L_k(\lambda) =  \lambda\max\left\{1,\ \frac{\nrm{\bbf}}{\nrm{\Gbf_{k}}}\right\}\\
U_k(\lambda) =  \frac{1}{\lambda}\min\left\{1,\ \frac{\nrm{\bbf}}{\nrm{\Gbf_{k}}}\right\}\end{dcases}, \qquad 1\leq k\leq IJ.
\end{equation}
The sparsity threshold $\widehat{\lambda}$ is then selected as the smallest minimizer of the cost function \begin{equation}\label{lossfcn}
\CalL(\lambda) = \frac{\nrm{\Gbf(\wbf(\lambda)-\wbf(0))}_2}{\nrm{\Gbf\wbf(0)}_2}+\frac{\nrm{\wbf(\lambda)}_0}{IJ}
\end{equation}
where $\wbf(\lambda):=\text{MSTLS}(\Gbf,\bbf; \lambda\,)$. We find $\widehat{\lambda}$ via grid search and set $\what=\text{MSTLS}(\Gbf,\bbf;\widehat{\lambda})$ as the output of the algorithm. In words, this is a modified sequential thresholding algorithm with non-uniform thresholds \eqref{MSTLSbnds} chosen based on the norms of the underlying library terms $\Gbf_{(i-1)J+j}\approx D^{\pmb{\alpha}^{(i)}}\psi*f_j(u)$ relative to the response vector $\bbf \approx D^{\pmb{\alpha}^{(0)}}\psi*u$. The purpose of this is to (a) incorporate relative sizes of library terms $\Gbf_{k}\wstar_k$ along with absolute sizes of coefficients $\wstar$ in the thresholding step, and (b) choose $\lambda$ automatically. 

\section{Opacity calculations}\label{sec:opacities}
Integrating the Planckian $B(\nu,T)$ (given in \eqref{eqB}) over photon frequency, we get
\begin{equation}\label{eq:B_int}
\int Bd\nu = \int_0^\infty \frac{2h\nu^3}{c^2}\left( e^{h\nu/kT} - 1\right)^{-1}d\nu = \frac{2 (kT)^4}{h^3c^2}\int_0^\infty x^3 (e^x-1)^{-1}dx = \frac{2\pi^4 (kT)^4}{15 h^3c^2}
\end{equation}
which serves to define the Stefan-Boltzmann constant $\sigma_\text{SB}$ and radiation constant $a$ through $\int Bd\nu = \frac{\sigma_\text{SB}}{\pi} T^4= \frac{ac}{4\pi} T^4$, or
\begin{equation}\label{eq:ac_relation}
ac := \frac{8\pi^5 k^4}{15 h^3c^2} = 4 \sigma_\text{SB}
\end{equation}
For the Larsen opacity 
\[\sigma(\nu,T) = \frac{\gamma}{(h\nu)^3}(1-e^{-h\nu/k T})\]
we have
\[\int \sigma Bd\nu = \int_0^\infty  \frac{\gamma}{(h\nu)^3}\frac{2h\nu^3}{c^2}\left( \frac{1-e^{-h\nu/k T}}{e^{h\nu/kT} - 1}\right)d\nu = \frac{2\gamma k T}{h^3c^2} \int_0^\infty e^{-x}dx = \frac{2\gamma k}{h^3c^2}  T\]
and thus
\begin{equation}\label{eq:larsen_sigp}
\sigma_P(T) = \frac{15\gamma}{\pi^4 k^3}T^{-3}.
\end{equation}
The rate of decrease in material temperature due to black-body radiation is thus linear in $T$: 
\begin{equation}\label{eq:larsen_sigp_term}
ac \sigma_p(T)T^4 = \frac{15 ac \gamma }{\pi^4 k^3}T = \frac{60 \sigma_\text{SB} \gamma }{\pi^4 k^3}T
\end{equation}

We define the quantity $\kappa_L$ in \eqref{eq:kappaL} using
\[\sigma_P(T_o;T_i) := \frac{\int \sigma(\nu,T_o)B(\nu,T_i) d\nu}{\int B(\nu,T_i)d\nu}\]
We compute 
\[\int \sigma(\nu,T_o)B(\nu,T_i) d\nu = \int_0^\infty  \frac{\gamma}{(h\nu)^3}\frac{2h\nu^3}{c^2}\left( \frac{1-e^{-h\nu/k T_o}}{e^{h\nu/kT_i} - 1}\right)d\nu = \frac{2\gamma k T_i}{h^3c^2} \int_0^\infty \left(\frac{1-e^{-\frac{T_i}{T_o}x}}{1-e^{-x}}\right) e^{-x}dx \]
Let $\eta = T_i/T_o$ and define $n=[\eta]$ to be the best integer approximation to $\eta$. We have that 
\[\int_0^\infty e^{-x}\left(\frac{1-e^{-\eta x}}{1-e^{-x}}\right) dx = \int_0^1 \left(\frac{1-y^\eta}{1-y}\right) dy \approx \int_0^1 \left(\frac{1-y^n}{1-y}\right) dy = \sum_{k=1}^n \frac{1}{k} \approx \log(n) + \gamma_\text{EM}\]
where $\gamma_\text{EM}$ is the Euler-Mascaroni constant. The approximation is very good so long as $T_i\gg T_o$. From this we have 
\[\sigma_P(T_o;T_i) \approx \frac{15}{\pi^4}\frac{\log(T_i/T_o)}{(k T_i)^3}\gamma\]

For the Rosseland weighted opacity, we have 
\[\int \partial_T B d\nu = \frac{2 h^2}{c^2k T^2}\int_0^\infty \nu^4e^{h\nu/k T} (e^{h\nu/k T}-1)^{-2} d\nu = \frac{2k^4T^3}{c^2 h^3}\int_0^\infty x^4 e^x(e^x-1)^{-2}dx = \frac{8\pi^4 k^4}{15 c^2 h^3}T^3\]
computed symbolically. The numerator is given by 
\[\int \frac{1}{\sigma}\partial_TB d\nu = \int_0^\infty\left(\frac{(h\nu)^3}{\gamma}(1-e^{-h\nu/k T})^{-1}\right)\left(\frac{2 h^2}{c^2k T^2} \nu^4e^{h\nu/k T} (e^{h\nu/k T}-1)^{-2} \right) d\nu\]
\[ = \frac{2h^5}{\gamma c^2 k T^2}\int_0^\infty \nu^7 e^{2h\nu/kT}(e^{h\nu/kT}-1)^{-3}d\nu = \frac{2k^7 T^6}{\gamma c^2 h^3} \int_0^\infty x^7 e^{2x}(e^x-1)^{-3} dx \approx (5.1047 \times 10^3)  \frac{2k^7}{\gamma c^2 h^3}T^6 \]
where the last integral is computed using \texttt{integral} in MATLAB. From this we get
\begin{equation}\label{eq:larsen_sigR}
\sigma_R(T) \approx \frac{4\pi^4}{15(5.1047 \times 10^3)} \frac{\gamma }{k^3}T^{-3} \approx  (5.0886\times 10^{-3})\frac{\gamma}{k^3}T^{-3}
\end{equation}

\bibliographystyle{ieeetr} %
\bibliography{refs.bib}




\end{document}